%% file: si2.tex
\documentclass[12pt,letterpaper]{amsart}
\usepackage{epic,eepic,latexsym, amssymb, amscd, amsfonts, xypic}

 \newlength{\baseunit}               
 \newcount{\numlines}                
\newcommand{\point}{\vspace{3mm}\par \noindent \refstepcounter{subsection}{\bf \thesubsection.} }
\newcommand{\tpoint}[1]{\vspace{3mm}\par \noindent \refstepcounter{subsection}{\bf \thesubsection.} 
  {\em #1. ---} }
\newcommand{\epoint}[1]{\vspace{3mm}\par \noindent \refstepcounter{subsection}{\bf \thesubsection.} 
  {\em #1.} }
\newcommand{\bpoint}[1]{\vspace{3mm}\par \noindent \refstepcounter{subsection}{\bf \thesubsection.} 
  {\bf #1.} }

\newcommand{\bpf}{\noindent {\em Proof.  }}
\newcommand{\epf}{\qed \vspace{+10pt}}

\newcommand{\R}{\mathbb{R}}
\newcommand{\C}{\mathbb{C}}

\newcommand{\G}{\mathbb{G}}

\newcommand{\proj}{\mathbb P}
\newcommand{\oh}{{\mathcal{O}}}

\newcommand{\al}{\alpha}

\newcommand{\be}{\beta}
\newcommand{\ga}{\gamma}

\newcommand{\De}{\Delta}
\newcommand{\si}{\sigma}

\newcommand{\Spec}{\operatorname{Spec}}

\newcommand{\init}{\rm{init}}
\newcommand{\final}{\rm{final}}

\newcommand{\bn}{\bullet_{\text{next}}}

\newcommand{\cb}{\circ \bullet}

\newcommand{\Xbb}{X_{\bullet} \cup X_{\bn}}

\newcommand{\cs}{\circ_{\text{swap}}}
\newcommand{\cn}{\circ_{\text{stay}}}

\newcommand{\cbs}{\cs \bn}
\newcommand{\cbn}{\cn \bn}
  
\newcommand{\Sv}{\Omega}

\newcommand{\bu}{\cdot}  

\newcommand{\Sch}{\mathbf{S}}
\newcommand{\bSv}{\pmb{\Omega}}
\newcommand{\Ggroup}{\rm{Gal}}
\newcommand{\Gvariety}{\rm{GalSch}}

\newcommand{\cited}{}

\newcommand{\secretnote}[1]{}

\newcommand{\lremind}[1]{{}}

\newcommand{\cut}[1]{}

\begin{document}
\pagestyle{plain}
\title{
{\Large {Schubert induction}}
}
\author{Ravi Vakil}
\address{Dept. of Mathematics, Stanford University, Stanford CA~94305--2125}
\email{vakil@math.stanford.edu}
\thanks{Partially supported by NSF Grant DMS--0228011, an AMS Centennial Fellowship, and an Alfred P. Sloan Research Fellowship.}
\date{Saturday, April 12, 2003.}
\begin{abstract}
We describe a {Schubert induction theorem}, a tool for analyzing
intersections on a Grassmannian over an arbitrary base ring.  The
key ingredient in the proof is the Geometric Littlewood-Richardson
rule of \cite{geolr}.

Schubert problems are among the most classical problems in enumerative
geometry of continuing interest.  As an application of Schubert
induction, we address several long-standing natural questions related
to {Schubert problems}, including: the ``reality'' of solutions;
effective numerical methods; solutions over algebraically closed
fields of positive characteristic; solutions over finite fields; a
generic smoothness (Kleiman-Bertini) theorem; and monodromy groups of
Schubert problems.  For example, we show that {\em all} Schubert
problems for {\em all} Grassmannians are enumerative over the real
numbers, completely answering the classical ``reality question'' for
Schubert problems.  These methods conjecturally extend to the flag
variety.
\end{abstract}
\maketitle
\tableofcontents

{\parskip=12pt 

\section{Introduction}

The Schubert induction theorem (Theorem~\ref{sit})
is a tool for studying intersections on a Grassmannian over
an arbitrary base, by deformations.
One motivation for such a result are {\em Schubert problems},
among the most classical problems in enumerative
geometry of continuing interest.  It is perhaps surprising
that many natural questions about Schubert problems
remain open.  We describe these questions, and give some answers, 
in Section~\ref{qa}.  
Most applications conjecturally extend to the flag variety 
(Sect.~\ref{flv}).
The theorem is stated and proved
in Section~\ref{pf}, and the  applications are  shown
in Sections~\ref{app} and~\ref{gm}.

\bpoint{Notation and philosophy} Fix a Grassmannian $G(k,n)= \G(k-1,n-1)$ over a base
field (or ring) $K$.  If $\al$ is a partition, let $\Sv_{\al} \in A^*(G(k,n))$ denote
the corresponding Schubert class.  Let
$\Sv_{\al}(F_{\bu})$  be the closed
 Schubert variety with respect to the flag $F_{\bu}$.  
Let $\bSv_{\al}(F_{\cdot}) 
\subset G(k,n) \times Fl(n)$ be the
universal Schubert variety.


Let $\pi_i: G(k,n) \times Fl(n)^m  \rightarrow G(k,n) \times Fl(n)$ ($1
\leq i \leq m)$ denote the projection, where the projection to $Fl(n)$ is
from the $i$th $Fl(n)$ of the domain.  We will make repeated use of the following diagram.
\lremind{starbucks}
\begin{equation}
\label{starbucks}
\xymatrix{
\pi_1^* \bSv_{\al_1}(F^1_{\bu})
\cap
\pi_2^* \bSv_{\al_2}(F^2_{\bu})
\cap \cdots
\cap
\pi_m^* \bSv_{\al_m}(F^m_{\bu})
\ar[d]^{{\ssize{\Sch}}} \ar[r] 
\ar@{^{(}->}[r] &
G(k,n) \times Fl(n)^m  \\
Fl(n)^m
}
\end{equation}

Questions about Schubert problems often reduce to questions about the
morphism $\Sch$, and in particular about $\Sch^{-1}(p)$ where $p$ is 
a general point of $Fl(n)^m$.  Schubert induction involves
studying this question by specializing $p$ through a carefully chosen
sequence of codimension 1 degenerations, 
ending with $p$ being a totally degenerate point
(where the $m$ flags coincide).  In each application, the property in
question will behave well under these degenerations.  The flavor of
the statement is that if something is true for $m=1$, then it is true
for all $m$.  This philosophy is described in detail
in Section~\ref{philosophy}.

As an informal but illustrative 
example, consider the statement ``there is a means of
combinatorially computing the number of preimages of $p$''.  Then
Schubert induction in this case amounts to the checker-tournament
method of solving Schubert problems described in\secretnote{See
  below.}  \cite[Sect.~3.12]{geolr}.  The base case is the trivial
statement ``there is one point in the zero-dimensional Schubert
variety''.

\bpoint{Acknowledgments} I am grateful to F.~Sottile for advice and
discussions.  His philosophy is clearly present in this paper.  In
particular, the phrase ``Schubert induction'' appeared first in
\cite[Sect.~1]{frankflag}, although with a slightly different meaning.
I thank A.~Buch and A.~Knutson for introducing me to this subject.  I
also thank B.~Sturmfels for suggesting that \cite{geolr} might have
numerical applications (see Sect.~\ref{qa} Question 5), D.~Allcock for
discussions on the symmetric group, B.~Poonen for advice on the
Chebotarev density theorem, and W.~Fulton for improving the exposition
of the main result.  I am grateful to H. Derksen for pointing out
(using the theory of quivers) that the Galois/monodromy group of
Schubert problems is sometimes not the full symmetric group, and for
producing explicit examples.

\section{Questions and answers} \label{qa} \lremind{qa}

Given a partition $\al$, the condition (i.e. element of $A^*(G(k,n))$ corresponding
to $\Sv_{\al}(F_{\bu})$ is called a {\em Schubert condition}.  A
{\em Schubert problem} is the following:

\bpoint{Schubert problem} Given $m$ Schubert conditions
$\Sv_{\al_i}(F^i_{\bu})$ with respect to general flags $F^i_{\bu}$ ($1
\leq i \leq m$) whose total codimension is $\dim G(k,n)$, what is the
cardinality of their intersection?  

In other words, how many $k$-planes satisfy various linear algebraic
conditions with respect to $m$ general flags?  (Or: what is the
cardinality of $\Sch^{-1}(F^1_{\bu}, \dots, F^m_{\bu})$ for general
$(F^1_{\bu}, \dots, F^m_{\bu}) \in Fl(n)^m$?)  This is the natural
generalization of the classical problem: how many lines in $\proj^3$
meet four (fixed) general lines?  The points of intersection are
called the {\em solutions of the Schubert problem}.  (For clarity's
sake, we say that the number of solutions is the {\em answer} to the
Schubert problem.)  An immediate (if imprecise) follow-up is: What can
one say about the solutions?

For example, if $K=\C$, the answer to the Schubert problems for $m=3$
are precisely the Littlewood-Richardson coefficients $c^{\ga}_{\al
  \be}$.

Suppose the base field is $K$,
and $\al_1$, \dots, $\al_m$ are given such that
$\dim \left( \Sv_{\al_1} \cup  \dots \cup \Sv_{\al_m} \right) = 0.$ 
The corresponding Schubert problem
is said to be {\em enumerative over $K$} if there are $m$ flags 
$F^1_{\bu}$, \dots, $F^m_{\bu}$ defined over $K$  such that
$\Sch^{-1}(F^1_{\bu}, \dots, F^m_{\bu})$ 
consists of 
$\deg \left( \Sv_{\al_1} \cup  \dots \cup \Sv_{\al_m} \right)$
(distinct) $K$-points.

\point \label{prototype} \lremind{prototype}
The answer to this problem over $\C$ is the prototype of the program
in enumerative geometry.  By the Kleiman-Bertini theorem
\cite{kleimanbertini}, the Schubert conditions intersect transversely,
i.e. at a finite number of reduced points.  Hence the problem is
reduced to one about the intersection theory of the Grassmannian.  The
intersection ring (the Schubert calculus) is known, using other
interpretations of the Littlewood-Richardson coefficients in
combinatorics or representation theory.

Yet many natural questions remain:

\bpoint{Reality questions} {} The classical ``reality question'' for
Schubert problems \cite[p.~55]{babyF}, \cite[Ch.~13]{it}, \cite[Sect.~9.8]{pr}
is:

\noindent {\bf Question 1.} 
{\em   Are all Schubert problems enumerative over $\R$?}

See \cite{santacruz, frankflag}
for this problem's history.  The case
$G(1,n)$ (and $G(n-1,n)$) is trivially linear algebra.  Sottile proved
the result for $G(2,n)$ (and $G(n-2,n)$) for all $n$, \cite{lines},
and for all problems involving only Pieri classes \cite{special}; see
\cite{jpaa} for further discussion.  The case
$G(2,n)$, as well as that of conics, also follows 
from \cite{crelle}.

This question can be fully answered with Schubert induction.

\tpoint{Proposition} {\em All Schubert problems for all Grassmannians are
  enumerative over $\R$.  Moreover, for a fixed $m$, there is a set of
  $m$ flags that works for {\rm all} choices of $\al_1$, \dots,
  $\al_m$.}
\label{enumerative} \lremind{enumerative}

This argument carries through with $\mathbb{R}$ replaced
by any field satisfying the implicit function theorem, such as
$\mathbb{Q}_p$.

As noted in \cite[Sect.~4.11(g)]{geolr}, Eisenbud's suggestion that
the deformations of the Geometric Littlewood-Richardson rule are a
degeneration of that arising from the osculating flag to a rational
normal curve, along with this proposition, would imply that the
Shapiro-Shapiro conjecture is true asymptotically.  (See \cite{eg} for
the proof in the case $k=2$.)

\bpoint{Enumerative geometry in positive characteristic} Enumerative
geometry in positive characteristic is almost a stillborn field,
because of the failure of the Kleiman-Bertini theorem.  (Examples of
the limits of our understanding are plane conics and cubics in
characteristic 2 \cite{vainsencher, berg}.)  In particular, the
Kleiman-Bertini Theorem fails in positive characteristic for all
$G(k,n)$ that are not projective spaces (i.e. $1<k<n-1$) --- Kleiman's
counterexample \cite[ex.~9]{kleimanbertini} for $G(2,4)$ easily
generalizes. D.~Laksov and R.~Speiser have developed a sophisticated
characteristic-free theory of transversality \cite{laksov, speiser,
  ls1, ls2}, but it does not apply in this case
\cite[Sect.~5]{sottilerecent}.

\noindent{\bf Question 2.} {\em 
Are Schubert problems enumerative over an algebraically closed field
of positive characteristic?}

To answer this question, we first answer a logically prior one:

\noindent{\bf Question 3.}   {\em
Is there any patch to the failure of the Kleiman-Bertini theorem
on Grassmannians?}

A related natural question is:

\noindent{\bf Question 4.} {\em 
Are Schubert problems enumerative over finite fields?}

We now answer all three questions.  The appropriate replacement of
Kleiman-Bertini is the following.  We say a morphism $f: X \rightarrow
Y$ is {\em generically smooth} if there is a dense open set $V$ of $Y$
and a dense open set $U$ of $f^{-1}(V)$ such that $f$ is smooth on
$U$.  If $X$ and $Y$ are varieties and $f$ is dominant, this is
equivalent to the condition that the function field of $X$ is
separably generated over the function field of $Y$.

\tpoint{Generic smoothness theorem}  \lremind{gst} \label{gst}
{\em The morphism $\Sch$ is generically smooth.
 More generally, if $Q
  \subset G(k,n)$ is a subvariety such that $Q \cap
  \bSv_{\al}(F_{\bu}) \rightarrow Fl(n)$ is generically smooth for
  all $\al$, then 
$$ \xymatrix{ Q \cap
\pi_1^* \bSv_{\al_1}(F^1_{\bu})
\cap
\pi_2^* \bSv_{\al_2}(F^2_{\bu})
\cap \cdots
\cap
\pi_m^* \bSv_{\al_m}(F^m_{\bu}) \ar[r] &  Fl(n)^m
}$$
 is as well.}

This begs the following question: is the only obstruction to the
Kleiman-Bertini theorem for $G(k,n)$ that suggested by Kleiman, i.e.
whether the variety in question intersects a general translate of all
Schubert varieties transversely?  More precisely, is it true that for
all $Q_1$ and $Q_2$ such that $Q_i \cap \bSv_{\al}(F_{\bu})
\rightarrow Fl(n)$ is generically smooth for all $\al$, and $i=1,2$,
it follows that
$$\xymatrix{
Q_1 \cap \sigma(Q_2) \ar[r] &  PGL(n)
}
$$
is also generically smooth, where $\sigma \in PGL(n)$?

Theorem~\ref{gst} answers Question 3, and 
leads to answers to Questions 2 and 4:

\tpoint{Corollary} {\em
\begin{enumerate} 
\item[(a)] All Schubert problems are enumerative for
  algebraically closed fields. 
\item[(b)] For any prime $p$, there is a
  positive density of points $P$ defined over finite fields of
  characteristic $p$ where $\Sch^{-1}(P)$ consists of $\deg \left(
    \Sv_{\al_1} \cup \dots \cup \Sv_{\al_m} \right)$ distinct
  points.
Moreover, for a fixed $m$, there is a positive density of 
points that works for {\em all} choices of $\al_1$ \dots, $\al_m$.
\end{enumerate}}

Part (a) follows as usual (see Sect.~\ref{prototype}). 
If $\dim \left(
  \Sv_{\al_1} \cup \dots \cup \Sv_{\al_m} \right) = 0,$ then
Theorem~\ref{gst} implies that $\Sch$ is generically separable
(i.e. the extension of function fields is separable).   Then (b)
follows by applying the Chebotarev density theorem for
function fields
to
$$ \xymatrix{
\coprod_{\al_1, \dots, \al_m} \pi_1^* \bSv_{\al_1}(F^1_{\bu})
\cap
\pi_2^* \bSv_{\al_2}(F^2_{\bu})
\cap \cdots
\cap
\pi_m^* \bSv_{\al_m}(F^m_{\bu})
\ar[r] &
Fl(n)^m
}$$
(see for example \cite[Lemma~1.2]{e}, although all that is needed
is the curve case, e.g. \cite[Sect.~5.4]{fj}).\secretnote{See Bjorn's
  e-mail of roughly Jan. 29, 2003.}

Sottile has proved transversality for intersection of
codimension 1 Schubert varieties \cite{sottilerecent}, and
P.~Belkale has recently proved transversality in general, 
using his proof of Horn's conjecture \cite[Thm.~0.9]{belkale}.

\bpoint{Effective numerical solutions (over $\C$) to all
   Schubert problems for all Grassmannians}
Even over the complex numbers, questions remain.

\noindent {\bf Question 5.}  {\em  Is there an effective numerical 
method for {\em solving} Schubert problems (i.e. calculating
the solutions to any desired accuracy)?}

 The case
of intersections of ``Pieri classes'' was dealt with in 
\cite{hss}. 
For motivation in control theory, see for example \cite{hv}.
In theory, one could numerically solve Schubert problems using the Plucker
embedding; however, this is unworkable in practice. 

Schubert induction leads to an
 algorithm for effectively numerically finding all solutions to all
 Schubert problems over $\C$.
The method will be described in \cite{sottilevakil}, and
the reasoning is sketched in Section~\ref{numsol}.

\bpoint{Galois or monodromy groups of Schubert problems} The Galois or
monodromy group of an enumerative problem measures three (related)
things:
\begin{enumerate} 
\item[(a)] (geometric) As the conditions are varied, how do the
solutions permute?
\item[(b)] (arithmetic) What is the field of definition of the solutions,
given the field of definition of the flags?
\item[(c)] (algebraic) What is the Galois group of the field extension
  of the ``variety of solutions'' over the ``variety of conditions''
  (see \eqref{starbucks})?
\end{enumerate}
The modern study of such problems was initiated by J. Harris in \cite{joe};
the connection between (a) and (c) is made there.
The connection to (b)
is via the Hilbert irreducibility theorem, as the
target of $\Sch$ is rational.

\noindent{\bf Question 6.} {\em 
What is the Galois group of a Schubert problem?}

We partially answer this question.  There is an explicit combinatorial
criterion that implies that a Schubert problem has Galois group ``at
least alternating'' (i.e. if there are $d$ solutions, the group is
$A_d$ or $S_d$).  This criterion holds over an arbitrary base ring.
To prove it, we will discuss useful methods for analyzing Galois
groups via degenerations.  The criterion is quite strong, and
seems to apply to all but a tiny proportion of Schubert problems.
For example:

\tpoint{Theorem}  {\em The Galois group of {\em any} Schubert problem
on the Grassmannians $G(2,n)$ ($n \leq 16$) and $G(3,n)$ ($n \leq 9$)
is either alternating
or symmetric.} \lremind{galeg} \label{galeg}

A short \verb+Maple+ program applying the criterion to a
general Schubert problem
is available upon request from the author.

One might reasonably expect that the Galois group of a Schubert
problem is always the full symmetric group.  However, this not the
case.  To our knowledge, the first examples are due to H.~Derksen.  In
Section~\ref{small} we describe the smallest example (involving 4
flags in $G(4,8)$), and determine (using the explicit checker
criterion) that the Galois action is that of $S_4$ on order $2$ subsets of $\{
1, 2,3,4 \}$.  In Section~\ref{family1} we give a family
of examples with $\binom N K$ solutions, with Galois group $S_N$, and
action corresponding to the $S_N$-action on order $K$ subsets of $\{ 1,
\dots, N \}$.

We also describe three-flag examples (i.e. corresponding to
Littlewood-Richardson coefficients) with similar behavior
(Sect.~\ref{family2}).  Littlewood-Richardson coefficients interpret
structure coefficients of the ring of symmetric functions as the
cardinality of some set.  These three-flag examples show that the set
has further structure, i.e.  the objects are not indistinguishable.
(More correctly, {\em pairs} of objects are not indistinguishable;
this corresponds to failure of two-transitivity.)

This family of examples was independently found by Derksen.  From his
quiver-theoretic point of view, the smallest member of this family (in
$G(6,12)$) corresponds to the extended Dynkin diagram of $E_6$, and
the smallest member of the other family (in $G(4,8)$) corresponds to
the extended Dynkin diagram of $D_4$.

\bpoint{Flag varieties}  \label{flv} \lremind{flv}
Conjecture~4.9\secretnote{conjecture}
of \cite{geolr} would imply that the results of this
paper except for those on Galois/monodromy groups apply to all
Schubert problems on flag manifolds.  In particular, as the conjecture
is verified for\secretnote{galveston} $n \leq 5$ 
\cite[Prop.~4.10]{geolr}, the results all hold in this range.  For
example:

\tpoint{Proposition} {\em All Schubert problems for $Fl(n)$ are
  enumerative over any algebraically closed field or any field with an
  implicit function theorem (e.g. $\R$) for $n \leq 5$.  For a fixed $m$,
  there is a set of $m$ flags that works for all choices of $\al_1$,
  \dots, $\al_m$.}

(The generalizations of the other statements in this paper are equally
straightforward.)

We note that in the case of triple intersections where the answer
is~1, Knutson has shown that the solution to the problem can be
obtained by using spans and intersections of the linear spaces in the
three flags \cite{kpc}; see also Purbhoo's result \cite{kevin}.

\section{The main theorem, and its proof}
\label{pf}

\bpoint{The key observation} 
\label{philosophy} \lremind{philosophy}
Suppose $f: X \rightarrow Y$ is a proper
morphism of irreducible varieties that we wish to show has some
property $P$.  We will require that $P$ satisfy several conditions,
including that it depend only on dense open subsets of the target
(condition {\bf (A)}).  An example of such a property is ``$f$ is
generically finite, and there is a Zariski-dense subset $U$ of real
points of $Y$ for which $f^{-1}(p)$ consists of $\deg f$ real points
for all $p \in U$.''

Suppose $D$ is a Cartier divisor of $Y$ such that $D \times_X Y$ is
reduced, and $D \times_X Y \rightarrow D$ has property $P$.  For good
choices of $P$ (call this condition on $P$ {\bf (C)}), such as the
example, this implies that $f$ has property $P$.

This motivates the following inductive approach.
Suppose
$$
X_0 = X \hookleftarrow X_1 \hookleftarrow X_2 \hookleftarrow \cdots
\hookleftarrow X_s$$
is a sequence of inclusions, where $X_{i+1}$ is a
Cartier divisor of $X_i$.  Suppose $Y_{i,j}$ ($1 \leq i \leq s$, $1
\leq j \leq J_i$) is a subvariety of $Y$ such that 
$f$ maps $Y_{i,j}$ to $X_i$, and $Y_{i,j}
\rightarrow X_i$ is proper, and for each $0 \leq i < s$, $1 \leq j \leq J_i$,
$$Y_{i,j} \times_{X_i} X_{i+1} = \cup_{j' \in I_{i,j}} Y_{i+1,j'}
$$
for some $I_{i,j} \subset J_{i+1}$, where each $Y_{i+1,j'}$ appears 
with multiplicity one.

If 
\begin{itemize}
\item $Y_{i+1,j'} \rightarrow X_{i+1}$ has $P$ for all $j' \in I_{i,j}$
implies $\cup_{j' \in I_{i,j}} Y_{i+1,j'} \rightarrow X_{i+1}$ has $P$
(condition {\bf (B)}), and 
\item $Y_{s,j} \rightarrow X_s$ has $P$ for all
$j \in J_s$ (condition {\bf (D)}, the base case for the induction),
\end{itemize}
then we may conclude that $f: X \rightarrow Y$ has $P$.  (Note that $Y
\times_X X_s \rightarrow X_s$ may be badly behaved; hence the need for
the inductive approach.)

The main result of this paper is that this process may be applied to
the morphism $\Sch$.  The key ingredient is the Geometric
Littlewood-Richardson rule.  

For some applications, we will need to refine the statement slightly.  For
example, to obtain lower bounds on monodromy groups, we will need the
fact that $I_{i,j}$ never has more than two elements.

\bpoint{Sketch of the Geometric Littlewood-Richardson rule \cite{geolr}}
The key ingredient in the proof of the Schubert induction theorem~\ref{sit}
is the Geometric Littlewood-Richardson rule, which we sketch here. 

The variety $Fl(n) \times Fl(n)$ is stratified by 
the locally closed subvarieties 
with fixed numerical data.  For each $(a_{ij})_{i,j \leq n}$, 
the corresponding subvariety is $\{
(F_{\bu}, F'_{\bu}) : \dim F_i \cap F'_j = a_{ij} \}$.
We denote such numerical data by the {\em configuration} $\bullet$ 
(normally interpreted as a
partition), and the corresponding locally closed subvariety
by $X_{\bullet}$.

The variety $ G(k,n) \times Fl(n) \times Fl(n) $ is the disjoint union
of ``two-flag Schubert varieties'', locally closed subvarieties with
specified numerical data.  For each $(a_{ij}, b_{ij})_{i,j \leq
  n}$, the corresponding subvariety is $\{ (F_{\bu}, F'_{\bu}, V) :
\dim F_i \cap {F'_j} = a_{ij}, \dim F_i \cap
F'_j \cap V = b_{ij}, \}$.  We denote the data of the
$(b_{ij})$ by $\circ$, so the locally closed subvarieties are indexed
by the configuration $\cb$.  Denote the corresponding two-flag
Schubert variety by $X_{\cb}$. (Warning: the closure of a two-flag
Schubert variety need not be a union of two-flag Schubert 
varieties\secretnote{cap} \cite[Cor.~ 3.13(a)]{geolr}, so this is not
a stratification.)

There is a {\em specialization order} $\bullet_{\init}$, \dots
$\bullet_{\final}$ in the Bruhat order, corresponding to partial
factorizations of the longest word\secretnote{seattle} \cite[Sect.~
2.2]{geolr}.  If $\bullet \neq \bullet_{\final}$ is in the
specialization order, then let $\bn$ be the next term in the order.
We have $X_{\bn} \subset \overline{X}_{\bullet}$, $\dim X_{\bn} = \dim
X_{\bullet} - 1$, $X_{\bullet_{\init}}$ is dense in $Fl(n) \times
Fl(n)$, and $X_{\bullet_{\final}}$ is the diagonal in $Fl(n) \times
Fl(n)$.

There is a subset of configurations $\cb$, called {\em mid-sort},
where $\bullet$ is in the specialization order.

\tpoint{Geometric Littlewood-Richardson Rule, inexplicit 
form\secretnote{checkergeometry} (cf. \cite[Sect.~3]{geolr})}
\lremind{inexplicit} \label{inexplicit} {\em 
\begin{enumerate}
\item[(i)] For any two partitions 
$\al_1$, $\al_2$, $\pi_1^* \bSv_{\al_1}(F^1_{\cdot})
\cap \pi_2^* \bSv_{\al_2} (F^2_{\cdot}) = \overline{X}_{\circ \bullet_{\init}}$ for some mid-sort
$\circ \bullet_{\init}$, or 
  $\al_2$, $\pi_1^* \bSv_{\al_1}(F^1_{\cdot})
\cap \pi_2^* \bSv_{\al_2} (F^2_{\cdot}) = \emptyset$.
\item[(ii)] For any mid-sort $\circ \bullet_{\final}$,
  $\overline{X}_{\circ \bullet_{\final}} = \pi_1^* \bSv_{\al} \cap
  \Delta \, (= \pi_2^* \bSv_{\al} \cap
  \Delta)$ for some $\al$, where $\Delta$ is the pullback to $G(k,n)
  \times Fl(n) \times Fl(n)$ of the diagonal $X_{\bullet_{\final}}$ 
 of $Fl(n) \times Fl(n)$.
\item[(iii)] For any mid-sort $\cb$ with $\bullet \neq
  \bullet_{\final}$, 
consider the diagram \cite[equ.~(1)]{geolr} \secretnote{bloop}\lremind{sibloop}
\begin{equation}
\label{sibloop}
\xymatrix{ 
\overline{X}_{\cb}  \ar@{^{(}->}[r]^{\text{open}} 
\ar[d] & 
\overline{X}_{\cb}  \ar[d] & D_X 
\ar@{_{(}->}[l]_{\text{closed}}
 \ar[d] 
\\
X_{\bullet}  \ar@{^{(}->}[r]^>>>>{\text{open}} &
X_{\bullet} \cup X_{\bn}  
& X_{\bn}.
\ar@{_{(}->}[l]_>>>>>{\text{closed}}
}
\end{equation}
The closures of $X_{\cb}$ are taken in $G(k,n) \times X_{\bullet}$
and $G(k,n) \times ( \Xbb )$ respectively, and 
the Cartier divisor $D_X$ is defined by fibered product. 
There are one or two mid-sort configurations (depending on $\cb$),
  denoted by $\cbs$ and/or $\cbn$, such that $D_X = 
  \overline{X}_{\cbs}$,
  $\overline{X}_{\cbn}$, or $\overline{X}_{\cbn} \cup \overline{X}_{\cbs}$.
\end{enumerate}}

There is a more precise version of this rule describing the mid-sort
$\cb$, and $\cbs$ and $\cbn$ (see\secretnote{checkergeometry}
\cite[Sect.~3]{geolr}).  For almost all applications here this version
will suffice, but the precise definition of mid-sort, $\cbs$, and
$\cbn$ will be implicitly required for the Galois/monodromy results
of Section~\ref{gm}.  

\bpoint{Statement of Main theorem}
Fix $Q \subset G(k,n)$, and define  $S = S(\al_1, \dots, \al_{m-1}) 
\subset G(k,n) \times Fl(n)^{m-1}$ by 
\begin{equation}
S := 
( Q \times Fl(n)^{m-1} ) \cap \pi_1^* \bSv_{\al_1}(F^1_{\bu}) \cap \cdots \cap
\pi_{m-1}^* \bSv_{\al_{m-1}}(F^{m-1}_{\bu}) .
\end{equation}
Then $S$ is irreducible, and the 
projection to $B := Fl(n)^{m-1}$
has relative dimension $\dim Q - \sum |\al_i|$.  (This follows easily
by constructing $S$ as a fibration over $Q$.)

Let $P$ be a property of morphisms depending only on dense
open subsets of the target, i.e. if $U \subset Y$ is a dense open
subset, then $f:  X \rightarrow Y$ has $P$ if and only if $f |_{f^{-1}(U)}$
has $P$ (call this condition {\bf (A)}).  For such $S \rightarrow B$,
and any mid-sort $\cb$, let $\rho_1$ and $\rho_2$ be the two
projections from $B \times (Fl(n)\times Fl(n))$ onto its factors.  Using
\eqref{sibloop}, construct \lremind{sibloop2}
\begin{equation}\label{sibloop2}
\xymatrix{
\rho_1^* S \cap \rho_2^* \overline{X}_{\cb}   
  \ar@{^{(}->}[r]^{\text{open}} 
\ar[d]^{\sssize (\dagger)} & 
\rho_1^* S \cap \rho_2^*  \overline{X}_{\cb}
 \ar[d] & 
\rho_1^* S \cap \rho_2^* D_X 
\ar@{_{(}->}[l]_{\text{closed}}
 \ar[d]_{\sssize (\dagger \dagger)} 
&
{\begin{array}{c} \sssize{\rho_1^* S \cap \rho_2^*  \overline{X}_{\cbs} }  \\
 \sssize{\text{ and/or }}  \\
\sssize{\rho_1^* S \cap \rho_2^*  \overline{X}_{\cbn}}  
\ar[l] \ar[dl]^{\sssize (\dagger \dagger \dagger)} \end{array}}
\\
B \times X_{\bullet}
  \ar@{^{(}->}[r]^>>>>{\text{open}} &
B \times \left( X_{\bullet} \cup X_{\bn} \right)
& B \times X_{\bn}.
\ar@{_{(}->}[l]_>>>>>{\text{closed}}
}
\end{equation}
As in \eqref{sibloop},
$\overline{X}_{\cb}$ is the closure of $X_{\cb}$ in the
appropriate space;\secretnote{ ($G(k,n) \times X_\bullet$ or $G(k,n) \times (\Xbb)$).}
$\rho_2^* \overline{X}_{\cb}$ is the pullback
of $\overline{X}_{\cb}$ from $X_{\bullet}$ or $\Xbb$, 
and similarly for the other terms of the top row.  
The upper right should be interpreted
as 
$$\rho_1^* S \cap \rho_2^* \overline{X}_{\cbs}, \quad  \rho_1^* S \cap \rho_2^* \overline{X}_{\cbn}, 
\quad \text{or} \quad
\rho_1^* S \cap \rho_2^* \overline{X}_{\cbs} \coprod \rho_1^* S \cap \rho_2^* 
\overline{X}_{\cbn},
$$
as in the Geometric 
Littlewood-Richardson rule~\ref{inexplicit}.

\tpoint{{\bf Schubert induction theorem}}  \label{sit} \lremind{sit}
{\em Suppose that for 
 any such $S \rightarrow B$ and any mid-sort $\cb$, 
{\bf (B)} if $(\dagger \dagger \dagger )$ has $P$, then $( \dagger \dagger) $ has $P$, and
{\bf (C)} if $(\dagger \dagger)$ has $P$, then $(\dagger)$ has $P$.
If the projection \lremind{sihyp}
\begin{equation}
\label{sihyp}
\xymatrix{ ( Q \times Fl(n))  \cap  \bSv_{\al}(F_{\bu}) \ar[r] &  Fl(n)}
\end{equation} 
has
$P$ for all partitions $\al$ (the ``base case'' of the Schubert
induction), then the projection 
$$
\xymatrix{ ( Q \times Fl(n)^m)  \cap \pi_1^* \bSv_{\al_1}(F^1_{\bu}) \cap \cdots \cap 
\pi_m^* \bSv_{\al_m}(F^m_{\bu}) \ar[r] &  Fl(n)^m}
$$
has $P$
for all $m$, $\al_1$, \dots, $\al_m$.}

In particular (taking $Q=G(k,n)$) if  the projection 
$$
\xymatrix{   \bSv_{\al}(F_{\bu}) \ar[r] &  Fl(n)}
$$
has $P$ 
for all $\al$ (condition {\bf (D)}), then the projection 
$$
\xymatrix{ \Sch:  \pi_1^* \bSv_{\al_1}(F^1_{\bu}) \cap \cdots \cap 
\pi_m^* \bSv_{\al_m}(F^m_{\bu}) \ar[r] & Fl(n)^m}
$$
has $P$.

\bpf  
We show that \lremind{toprove}
\begin{equation}
\label{toprove}
( Q \times Fl(n)^{m-1} \times X_{\bullet}) \cap \pi_1^* \bSv_{\al_1}(F^1_{\bu}) \cap \cdots \cap
\pi_{m-1}^* \bSv_{\al_{m-1}}(F^{m-1}_{\bu}) \cap \rho^* \overline{X}_{\cb} 
\rightarrow
Fl(n)^{m-1} \times X_{\bullet}
\end{equation}
(where $\rho$ is the projection to $X_{\bullet}$) 
has $P$ for all $m$ and mid-sort
$\cb$, by induction on $(m, \bullet)$, where $(m_1, \bullet_1)$
precedes $(m_2, \bullet_2)$ if $m_1 < m_2$, or $m_1=m_2$ and
$\bullet_1 < \bullet_2$ in the specialization order.  

{\em Base case $m=1$, $\bullet = \bullet_{\final}$.}
By the Geometric Littlewood-Richardson rule~\ref{inexplicit} (ii),
$$
\xymatrix{
( Q \times X_{\bullet_{\final}})  
\cap \rho^* \overline{X}_{\circ \bullet_{\final}} 
\ar[r] \ar@{<->}[d]_{ {\sssize \cong}} &
X_{\bullet_{\final}} \ar@{=}[d]
&  \text{(i.e., \eqref{toprove})} \\
( Q \times X_{\bullet_{\final}})  \cap \pi_1^* \bSv_{\al} \cap \Delta 
\ar[r] \ar@{<->}[d]_{{\sssize \cong}} &
 X_{\bullet_{\final}} 
\ar@{<->}[d]^{{\sssize \cong}}  \\
( Q \times Fl(n))  \cap  \bSv_{\al}(F_{\bu})  \ar[r] & 
Fl(n)
}
$$
has $P$ by \eqref{sihyp}.

{\em Inductive step, case $\bullet \neq \bullet_{\final}$.}
By the inductive hypothesis,
$$
( Q \times Fl(n)^{m-1} \times X_{\bn}) \cap \pi_1^* \bSv_{\al_1}(F^1_{\bu}) \cap \cdots \cap
\pi_{m-1}^* \bSv_{\al_{m-1}}(F^{m-1}_{\bu}) \cap \rho^* \overline{X}_{\cbn} 
\rightarrow
Fl(n)^{m-1} \times X_{\bn}$$
and/or
$$
( Q \times Fl(n)^{m-1} \times X_{\bn}) \cap \pi_1^* \bSv_{\al_1}(F^1_{\bu}) \cap \cdots \cap
\pi_{m-1}^* \bSv_{\al_{m-1}}(F^{m-1}_{\bu}) \cap \rho^* \overline{X}_{\cbs} 
\rightarrow
Fl(n)^{m-1} \times X_{\bn}$$
have $P$.  
Then an application of {\bf (B)} and {\bf (C)} shows that \eqref{toprove}
has $P$ as well.

{\em Inductive step, case $\bullet = \bullet_{\final}$, $m>1$.}
Suppose $\overline{X}_{\circ \bullet_{\final}} = \pi_1^* \bSv_{\al} \cap \Delta$
and  $\pi_1^* \bSv_{\al_{m-1}}(F_{\cdot})
\cap \pi_2^* \bSv_{\al} (F'_{\cdot}) = \overline{X}_{\circ' \bullet_{\init}}$
(using the Geometric Littlewood-Richardson rule~\ref{inexplicit} (ii) and (i)
respectively).  Then
$$
\xymatrix{
\scriptstyle{  ( Q \times Fl(n)^{m-1} \times X_{\bullet}) \cap \pi_1^* \bSv_{\al_1}(F^1_{\bu}) \cap \cdots \cap
 \pi_{m-1}^* \bSv_{\al_{m-1}}(F^{m-1}_{\bu})
 \cap 
 \rho^* \overline{X}_{\cb} 
}
\ar[r] 
\ar@{<->}[d]_{ {\sssize \cong}} &
\scriptstyle{
Fl(n)^{m-1} \times X_{\bullet}  
}
\ar@{<->}[d]^{{\sssize \cong}}  
\\
\scriptstyle{
( Q \times Fl(n)^{m-1} \times X_{\bullet}) \cap \pi_1^* \bSv_{\al_1}(F^1_{\bu}) \cap \cdots \cap
\pi_{m-1}^* \bSv_{\al_{m-1}}(F^{m-1}_{\bu}) 
\cap  \pi_m^* \bSv_{\al}(F^m_{\bu})
}
\ar[r] \ar@{<->}[d]_{{\sssize \cong}} &
\scriptstyle{
Fl(n)^m
}
 \ar@{=}[d]
\\
\scriptstyle{
( Q \times Fl(n)^{m-1} \times X_{\bullet}) \cap \pi_1^* \bSv_{\al_1}(F^1_{\bu}) \cap \cdots \cap
\pi_{m-2}^* \bSv_{\al_{m-2}}(F^{m-2}_{\bu})
 \cap \rho^* \overline{X}_{\circ' \bullet_{\init}} 
}
 \ar[r] & 
\scriptstyle{
Fl(n)^m
}
}
$$
which has $P$ as  (by {\bf (A)})
$$
( Q \times Fl(n)^{m-1} \times X_{\bullet}) \cap \pi_1^* \bSv_{\al_1}(F^1_{\bu}) \cap \cdots \cap
\pi_{m-2}^* \bSv_{\al_{m-2}}(F^{m-2}_{\bu}) \cap \rho^* {\overline{X}}_{\circ' \bullet_{\init}} 
 \rightarrow  
 Fl(n)^{m-2} \times X_{\bullet_{\init}}
$$
has $P$ by the inductive hypothesis.\secretnote{In that last line, the overline on 
${X}_{\circ' \bullet_{\init}}$ is a closure in a different space.} 
\epf

For some applications, we will need a slight variation.

\tpoint{{\bf Schubert induction theorem, bis}} \lremind{sit2}
\label{sit2}
{\em Suppose $P$ satisfies conditions {\bf (A--C)}.  
If
$$  \xymatrix{
\coprod_{\al: \; \dim Q - |\al|=0} ( Q \times Fl(n))  \cap  \bSv_{\al}(F_{\bu}) \ar[r] & Fl(n)
}$$
has
$P$ then 
$$ \xymatrix{
\coprod_{\al_1, \dots, \al_m : \; \dim Q - \sum |\al_i|=0}
( Q \times Fl(n)^m)  \cap \pi_1^* \bSv_{\al_1}(F^1_{\bu}) \cap \cdots \cap 
\pi_m^* \bSv_{\al_m}(F^m_{\bu}) \ar[r] & Fl(n)^m
}$$
 has $P$
for all $m$.}

The proof is identical to that of Theorem~\ref{sit}.

\section{Applications} \label{app} 
\lremind{app}

We now verify the conditions {\bf (A--C)} for several $P$
to prove the results claimed in Section~\ref{qa}.

\bpoint{Positive characteristic:
Proof of Proposition~\ref{gst}}
Let $P$ be the property that the morphism $f$ is generically smooth.
Then $P$ clearly satisfies {\bf (A--C)} (note that the relative
dimensions of ($\dagger$ -- $\dagger \dagger \dagger$) are the same,
and that $X_{\cbn}$ and $X_{\cbs}$ are disjoint),
and the Schubert induction hypothesis {\bf (D)}; apply
Theorem~\ref{sit}.  \epf


\bpoint{Reality: Proof of Proposition~\ref{enumerative}} Let $P$ be
the property that there is a Zariski-dense subset $U$ of real points
of $Y$ for which $\Sch^{-1}(p)$ consists of $\deg \Sch$ real points
for all $p \in U$.  Clearly $P$ satisfies {\bf (A--C)} and the
Schubert induction hypothesis {\bf (D)}.  Apply Theorem~\ref{sit2}.
\epf

As mentioned earlier, the same argument applies to any field
satisfying the implicit function theorem, such as $\mathbb{Q}_p$.

\bpoint{Numerical solutions} \label{numsol} \lremind{numsol}
Informally, this  application corresponds to applying Theorem~\ref{sit2}
to the property
of generically finite morphisms $f: X \hookrightarrow G(k,n) \times Y
\rightarrow Y$ ``for each point of $Y$ whose preimage is a finite number of
points, there is an effective algorithm for numerically finding these
points''.  Condition {\bf (C)} corresponds to the fact that if
$|f^{-1} (y)| = \deg f$ and the points of $f^{-1}(y)$ 
can be numerically calculated, then by the implicit function theorem,
the points of $f^{-1}(y')$ can be numerically calculated for 
all $y'$ such that $\dim f^{-1}(y') = 0$.  This
idea will be developed in \cite{sottilevakil}.

\section{Galois/monodromy groups of Schubert problems}
\label{gm}
\lremind{gm}

\point \label{tree}\lremind{tree}
We recall the ``checker tournament'' 
algorithm\secretnote{not sure; ``Enumerative problems and
checker tournaments}
\cite[Sect.~3.12]{geolr}
for solving Schubert problems.
We begin with $m$ partitions, and we make a series
of moves.  Each move consists of one of the following.
\begin{enumerate}
\item[(i)] Take two partitions,  and begin a checker game if possible,
else end the tournament.
\item[(ii)] Translate a completed checker game back to a partition.
\item[(iii)]  Make a move in an ongoing checkergame.
\end{enumerate}
These parallel (i)--(iii) of Theorem~\ref{inexplicit}.
When one partition and no checker games are left, the 
tournament is complete.  At step (iii), the
checker tournament may bifurcate (if both a ``stay''
and a ``swap'' are possible), and both branches must
be completed.

This answer to the Schubert problem can be interpreted as a
creating a directed tree, where the vertices correspond to partially
completed checker tournament.  Each vertex has in-degree 1 (one
immediate ancestor) except for the root (corresponding to the original
Schubert
problem), and out-degree (number of immediate descendants) between 0
and 2.  The graph is
constructed starting with the root, and for each vertex that is not a
completed checkergame, a choice may be made (depending on (iii)) which
may lead to a bifurcation.  
Vertices corresponding to a single partition and no checkergames
are called {\em leaves}.  (There may be other vertices with out-degree 0,
arising from (i); these are not leaves.)
The answer is the number of
leaves of the tree.

The answer is of course independent of the choices made; in the
description of\secretnote{not sure, same as previous}
\cite[Sect.~3.12]{geolr}, and in the proof of Theorems~\ref{sit} and
\ref{sit2}, each checker game was chosen to be completed before the
next was begun.  

\tpoint{Theorem}  \lremind{monthm} \label{monthm}
{\em Suppose we are given a Schubert problem 
such that there is a directed tree as above, where
each vertex with out-degree two  satisfies either
\begin{enumerate}
\item[(a)]  there are a different number of leaves on
the 2 branches, or
\item[(b)]  there is one leaf on each branch.
\end{enumerate}
Then the Galois group of the Schubert problem is at least 
alternating.}

\bpoint{Specialization of monodromy} To prove the theorem, we will
examine how Galois groups behave under specialization.  

We say a generically finite morphism $f: X \rightarrow Y$ is {\em
  generically separable} if the corresponding extension of function
fields is separable.  Define the Galois group $\Ggroup_f$ of a
generically finite and separable (i.e. generically \'{e}tale)  morphism to
be the Galois group of the Galois closure of the corresponding
extension of function fields.

\epoint{Remark:  The complex case}\label{complexcase} \lremind{complexcase}
To motivate  later statements over an arbitrary ground ring, we
first consider the complex case.  Suppose
$$
\xymatrix{
W \ar[d] & Y \ar@{_{(}->}[l]_>>>>>{\text{closed}} \ar[d] \\
X        & Z \ar@{_{(}->}[l]_>>>>>{\text{closed}} }
$$
is a fiber diagram of complex schemes, where the vertical morphisms
are proper generically finite degree $d$; $W$, $X$, and $Z$ are
irreducible varieties; $Z$ is Cartier in $X$; $X$ is regular in
codimension 1 along $Z$; and $Y$ is reduced.  Then $\Ggroup_{W
  \rightarrow X}$ can be interpreted as an element of $S_d$ by fixing
a point of $X$ with $d$ preimages, and considering loops in the smooth
locus of $X$ based at that point, and their induced permutations of
the preimages.

(a) If $Y$ is irreducible, then by interpreting $\Ggroup_{Y \rightarrow
  Z}$ by choosing a general base point of $Z$ and elements of the
fundamental group of the smooth part of $Z$ generating the Galois
group, 
we have
constructed
an inclusion $\Ggroup_{Y \rightarrow Z} \hookrightarrow \Ggroup_{W
  \rightarrow X}$.  In particular, if the first group is at least
alternating, then so is the second.

(b) If $Y$ has two components $Y_1$ and $Y_2$, which each map generically
finitely onto $Z$ with degrees $d_1$ and $d_2$ respectively (so
$d_1+d_2=d$), then the same construction produces a subgroup $H$ of
$\Ggroup_{Y_1 \rightarrow Z} \times \Ggroup_{Y_2 \rightarrow Z}$ which
surjects onto $\Ggroup_{Y_i \rightarrow Z}$ (for $i=1,2$), and an
injection of $H$ into $\Ggroup_{W \rightarrow X}$ (via the
induced 
inclusion $S_{d_1} \times S_{d_2} \hookrightarrow S_d$).

Then  a purely group-theoretical argument 
(Prop.~\ref{groupie}) relying on
Goursat's lemma will show that if $\Ggroup_{Y_i \rightarrow Z}$ is at
least alternating ($i=1,2$), $W$ is connected (so $\Ggroup_{W \rightarrow
X}$ is transitive), and $d_1 \neq d_2$ or $d_1=d_2=1$, then $\Ggroup_{W \rightarrow X}$
is at least alternating as well.

\epoint{The general case} With this complex intuition in hand, we
prove Remarks~\ref{complexcase}(a) and (b) over an arbitrary ring.
Suppose $k_1 \subset k_2$ is a separable degree $d$ field extension.
Choose an ordering $x_1$, \dots, $x_d$ of the $\overline{k}_1$-valued
points (over $\Spec k_1$) of $\Spec k_2$.  If $g: X \rightarrow Y$ is
a generically finite separable (i.e. generically \'{e}tale) morphism,
define the ``Galois scheme'' $\Gvariety_g$ by $$\overbrace{X \times_Y
  \cdots \times_Y X }^{\deg g} \setminus \De$$
where $\De$ is the
``big diagonal''.  Recall that the Galois group of $k_1 \subset k_2 $
can be interpreted as a subgroup of $S_d$ as follows: $\si$ is in the
Galois group if and only if $(x_{\si(1)}, \dots, x_{\si(d)})$ is in
the same component of $\Gvariety_{\Spec k_2 \rightarrow \Spec k_1}$ as
$(x_1, \dots, x_d)$.

To understand how this behaves in families, let $R$ be a discrete
valuation ring with function field $K$ and residue field $k$.  Suppose
the following is a fiber diagram
$$
\xymatrix{
X_K \ar[d] 
  \ar@{^{(}->}[r]^>>>>>{\text{open}} 
&  X_R \ar[d]  &  X_k \ar[d] 
\ar@{_{(}->}[l]_>>>>>>{\text{closed}}   \\
\Spec K 
  \ar@{^{(}->}[r]^>>>>{\text{open}} 
&  \Spec R &  \Spec k 
\ar@{_{(}->}[l]_>>>>>{\text{closed}}  \\
}
$$
where $X_K$ is irreducible, $X_k$ is reduced, and the vertical morphisms
are finite and separable (and hence \'{e}tale, using $X_k$ reduced).  

After choice of algebraic closures, there is a bijection from the
$\overline{K}$-valued points of $X_K$ with the $\overline{k}$-valued
points of $X_k$.

By observing that each component of $\Gvariety_{X_k \rightarrow \Spec
  k}$ lies in a unique irreducible component of $\Gvariety_{X_R
  \rightarrow \Spec R}$ (as $\Gvariety_{X_R \rightarrow \Spec R}
\rightarrow \Spec R$ is \'{e}tale along $X_k$), we see that Remarks~\ref{complexcase}
(a) and (b) hold in general (by applying these comments to $\oh_{X,Z}$).

\epoint{Proof of Theorem~\ref{monthm}}
Each vertex $v$ corresponds to a diagram
$$
\xymatrix{
X  \ar@{^{(}->}[r] \ar[d]_{f_v} & B \times G(k,n) \\
B
}
$$
where $X$ is irreducible, and $B$ is a product of flag varieties
(one for each partition) and strata $X_{\cb}$ (one for each
checkergame in process).  The morphism $f_v$ is generically finite and
separable, and its degree is the answer
to the corresponding enumerative problem.

We label each vertex $v$ with the number of leaves on that branch
(i.e.  with $\deg f_v$), and with the Galois group $\Ggroup_v$ of that
problem.  We prove that $\Ggroup_v = A_{\deg f_v}$ or $S_{\deg f_v}$
for all $v$ by induction on $v$.  (This is a slight generalization of
Schubert induction.)  

If $v$ is a leaf, the result is trivial.

Suppose next that $v$ is a vertex with one descendant $w$, and $\Ggroup_w$
is at least alternating.  If the move from $v$ to $w$ is of type (i)
or (ii), then the morphism $f_v$ is the same as $f_w$, so the result holds.  If the move from vertex $v$ is
of type (iii) (so exactly one of $\{ \text{stay}, \text{swap} \}$ is
possible), then $G_v$ is at least alternating by
Remark~\ref{complexcase}(a).

Next, suppose $v$ has two immediate descendants, so we are in case
(iii) and both ``stay'' and ``swap'' are possible. 
If one branch has no leaves, then $G_v$ is at least alternating
by Remark~\ref{complexcase}(a), so assume otherwise.
As $X$ is irreducible, $\Ggroup_v$ is transitive.
By Remark~\ref{complexcase}(b) and the group theoretic calculation
of Proposition~\ref{groupie}, $\Ggroup_v$ is at least alternating, and the
inductive step is complete.

Thus by induction the root vertex has at least alternating Galois
group, completing the proof of the Theorem~\ref{monthm}.\epf

\tpoint{Proposition}  \label{groupie} \lremind{groupie}
{\em Suppose $G$ is a transitive subgroup of $S_{m+n}$ 
such that $G \cap (S_m \times S_n)$ contains
a subgroup $H$ such that the projection of $H$ to
$S_m$ (resp. $S_n$) is either $A_m$ ($m \geq 4$) 
or $S_m$ (resp.
$A_n$ for $n \geq 4$,  or $S_n$).  
\begin{enumerate}
\item[(a)]  
If $m \neq n$, then $G= A_{m+n}$ ($m+n \geq 4$) or $S_{m+n}$.
\item[(b)] If $m=n=1$ then $G=S_2$.
\end{enumerate}
}

Note that if $m=n$, then 
$$
\{ e, (1,n+1)(2,n+2) \cdots (n,2n) \} \rtimes 
(S_{\{1, \dots, n \}} \times
S_{ \{ n+1, \dots, 2n  \}})$$
is a subgroup of $S_{2n}$ whose intersection with
$S_{\{1, \dots, n \}} \times
S_{ \{ n+1, \dots, 2n  \}}$ surjects onto
each of its factors.

\bpf  Part (b) is trivial, so we prove (a).
Assume without loss of generality that $n>m$.

  Recall Goursat's lemma:\secretnote{Lang 2nd ed. p. 54 Ex. 9} 
if $H \subset G_1 \times G_2$,
such that $H$ surjects onto both factors, then there are normal
subgroups
$N_i \lhd G_i$ ($i=1,2$) and an isomorphism
$\phi: G_1/N_1 \stackrel \sim \rightarrow G_2/N_2$ such that
$(g_1, g_2) \in H$ if and only if $\phi(g_1N_1) =g_2 N_2$.

We first show that if $G$ is a transitive subgroup of $S_{m+n}$
($n> m \geq 3$) containing $A_m \times A_n$, then $G$ contains any
3-cycle and hence $A_{m+n}$.  Color the numbers $1$ through $m$ red
and $m+1$ through $m+n$ green.  Any monochromatic 3-cycle lies in $A_m
\times A_n$ and hence $G$.  Suppose $\tau$ is any element of $G$
sending a green number to a red position.  (i) If there are two
numbers of each color in the positions of one color (say, red) then
the conjugate of a 3-cycle in $A_m$ by $\tau$ will be a 3-cycle $\al$
of 1 red and 2 green objects, and the conjugate of a different
3-cycle in $A_m$ by $\tau$ will be a 3-cycle $\al$ of 1 green and 2
red objects.  Similarly, (ii) if there is at least one number of each
color in the positions of {\em both} colors, we can find a conjugate
$\al$ of a 3-cycle in $A_m$ or $A_n$ by $\tau$ that is a 3-cycle $\al$
of 1 red and 2 green objects, and the conjugate of a 3-cycle in $A_n$
or $A_m$ by $\tau$ that is a 3-cycle $\al$ of 1 green and 2 red
objects.  By conjugating $\al$ and $\be$ further by elements of $A_m
\times A_n$, we can obtain any non-monochromatic 3-cycle.  Now $\tau$
falls into case (i) and/or (ii), or $n=m+1$ and $\tau$ sends all red
objects to green positions, and all but one green object to red
positions.  Suppose $p$ is the green position containing the green
object in $\tau$, and $\si$ is a 3-cycle in $A_n$ moving $p$.  Then
$\tau^{-1} \si \tau$ is a permutation where exactly one red object is
sent to a green position, and vice versa, and we are in case (ii). Thus
in all cases $G$ contains $A_{m+n}$, as desired.

We now deal with the case $m,n \geq 3$.  By Goursat's lemma, $G$ must
contain $A_m \times A_n$.  (For example, if the projections of $H$ to
$S_m$ and $S_n$ are surjective, then $H$ arises from isomorphic
quotients $S_n/N_m \cong S_n / N_n$.  Then $(N_m,N_n) = (S_m,S_n)$ or
$(N_m,N_n)= (A_m,A_n)$; in both cases $A_m \subset N_m$ and $A_n
\subset N_n$.)  Then apply the previous paragraph.

For the remaining cases, it is straightforward to see (using Goursat)
that (i) if the image of $H$ is $A_n$ (resp. $S_n$) and $m=1$, then
$A_{m+n} \subset G$ (resp. $G=S_{m+n}$), and (ii) if $H$ surjects onto
$A_n$ ($n \geq 4$) or $S_n$ and $m=2$, then $A_{m+n} \subset G$. \epf

\epoint{Remark} \label{smallrk} \lremind{smallrk} We note for use in
Section~\ref{small} that if $n=1$ and the projection to $S_m$ is
surjective, then the same argument shows that $G=S_{m+n}$.

\bpoint{Applying Theorem~\ref{monthm}} Theorem~\ref{monthm} is quite
strong, and can be checked with a naive  computer program.  For
example, it implies that all Schubert problems for $G(2,n)$ for $n
\leq 16$ are at least alternating.  It also
implies that all but a tiny handful of Schubert problems for Grassmannians
of dimension less than 20 are at least alternating; we will describe
these exceptions.  

For $k>1$, the criterion will fail for the Schubert problem
$(1)^{k^2}$ on $G(k,2k)$: the first degeneration (i.e. the first
vertex with out-degree 2) will correspond to
$$
(1)^{k^2} = (2) (1)^{k^2} + (1,1) (1)^{k^2}
$$
and the two branches will have the same number of leaves by symmetry.
More generally, if $1 \leq m<k$ and $(m,k) \neq (1,2)$, 
the criterion will fail for the Schubert problem
$$
(\overbrace{m, \dots, m}^m) (1)^{k^2-m^2}$$
on $G(k,2k)$ for the same reason.

\begin{figure}[ht]
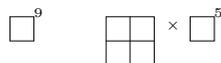

\begin{center}
\include{g36}
\end{center}
\caption{The two counterexamples of $G(3,6)$}
\label{g36}
\end{figure}

\begin{figure}[ht]
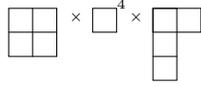

\begin{center}
\include{g37a}
\end{center}
\caption{An induced (non-primitive) counterexample in $G(3,7)$}
\label{g37a}
\end{figure}

On $G(3,6)$, the only counterexamples are of this sort, when $m=1$ and
$2$, shown in Figure~\ref{g36}.  By ``embedding'' these problems in
larger problems, these trivially induce counterexamples in larger
Grassmannians; for example, Figure~\ref{g37a} is a counterexample in
$G(3,7)$ that is really an avatar of the second example in $G(3,6)$.
We call counterexamples in $G(k,n)$ not arising in this way, i.e.
involving only subpartitions not meeting the right column and bottom
row of the rectangle, {\em primitive} counterexamples.

\begin{figure}[ht]
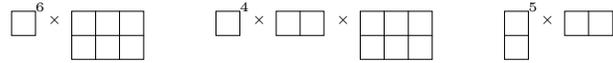

\begin{center}
\include{g37}
\end{center}
\caption{The primitive counterexamples in $G(3,7)$}
\label{g37}
\end{figure}

Then $G(3,7)$ has only three counterexamples, shown in
Figure~\ref{g37}, and the counterexamples in $G(4,7)$ are given by the
transposes of these.  The Grassmannian $G(3,8)$ has six
counterexamples,\secretnote{See jan2103.tex} shown in
Figure~\ref{g38}, and $G(3,9)$ has 13 counterexamples,\secretnote{See
  handwritten notes Jan. 25} shown in Figure~\ref{g39}.

\begin{figure}[ht]
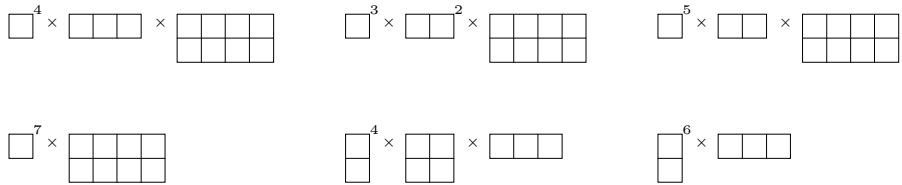

\begin{center}
\include{g38}
\end{center}
\caption{The primitive counterexamples in $G(3,8)$}
\label{g38}
\end{figure}

\begin{figure}[ht]
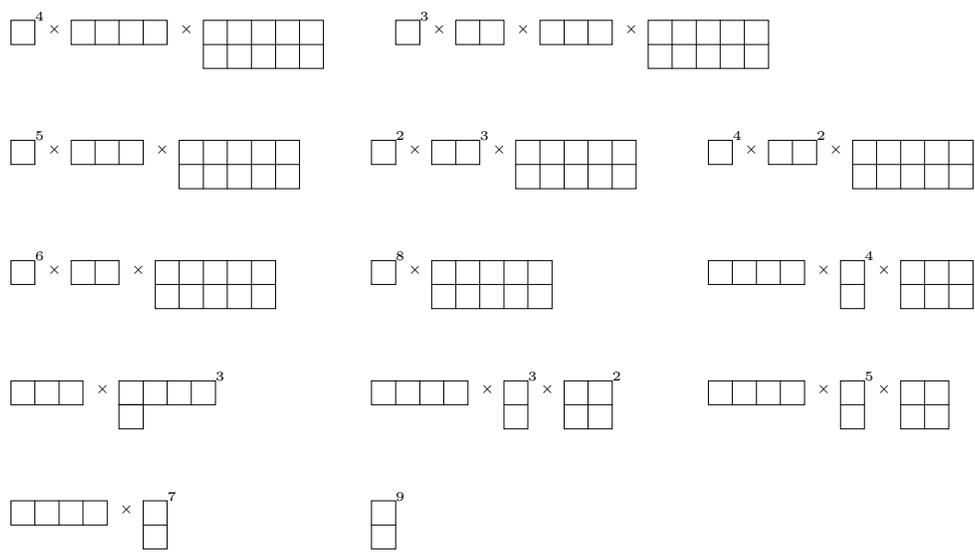

\begin{center}
\include{g39}
\end{center}
\caption{The primitive counterexamples in $G(3,9)$}
\label{g39}
\end{figure}

All of these exceptions can be excluded with the following,
slightly stronger criterion.

\tpoint{Theorem}  \lremind{monthm2} \label{monthm2}
{\em Suppose we are given a Schubert problem 
such that there is a directed tree as above, where
each vertex with two immediate descendants satisfies (a)
or (b) of Theorem~\ref{monthm}, or 
\begin{enumerate}
\item[(c)] there are $m \neq 6$ leaves on each branch,
and it is known that the corresponding Galois group
is two-transitive.
\end{enumerate}
Then the Galois group of the Schubert problem is at least 
alternating.}

In particular, to show that the Galois group is $(n-2)$-transitive,
it often suffices to show that it is two-transitive.

As with Theorem~\ref{monthm}, the proof reduces to the following
variation of Proposition~\ref{groupie}.

\tpoint{Proposition} \label{groupie2} \lremind{groupie2} {\em Suppose
  $G$ is a two-transitive subgroup of $S_{2m}$ ($m \neq 6$) such that
  $G \cap (S_m \times S_m)$ contains a subgroup $H$ such that the
  projection of $H$ to both factors $S_m$ is either $A_m$ ($m \geq 4$)
  or $S_m$.  Then $G=A_{2m}$ or $S_{2m}$.  }

The proof is similar to that of Proposition~\ref{groupie}, and
is omitted.

If $m=n=6$, D.~Allcock has pointed out that the Mathieu group
$M_{12}$ can be expressed as a subgroup of $S_{12}$ such that
$$M_{12} \cap (S_6 \times S_6) = \{ (g, \si(g) ) : g \in S_6 \}$$
where $\si$ is an outer automorphism of $S_6$.  Thus Proposition~\ref{groupie2} cannot be extended to $m=6$.

We say two vertices $v$, $w$ in a directed tree (as in
Sect.~\ref{tree}) are {\em equivalent} if they are connected by
a chain of edges $v_1 \rightarrow v_2 \rightarrow \cdots \rightarrow
v_s$ ($(v_1, v_s) = (v,w)$ or $(w,v)$) and $\deg f_v = \deg f_w$ (and hence
$= \deg f_{v_i}$ for all $i$).  
In each of the cases $G(3,n)$ ($6 \leq n \leq 9$)
given above, it is possible to find such a tree satisfying
Theorem~\ref{monthm2} (a)--(c), where the vertices of type (c) are
equivalent to vertices corresponding to Schubert problems (i.e.
corresponding to a set of partitions, with no
checkergames-in-progress), and to show by ad hoc means that these
Schubert problems are two-transitive.  (The details are omitted; this
method should not be expected to be workable in general.)  Hence all
Schubert problems for these Grassmannians have Galois group at least
alternating.

The Grassmannian $G(4,8)$ has only 31 Schubert problems where the
criterion of Theorem~\ref{monthm} does not apply (not shown here).
Each of these cases may be reduced to checking that a certain Schubert
problem is two-transitive.  As we shall see in the next Section,
in one of these cases two-transitivity does not hold!

\bpoint{Galois groups of Schubert problems needn't be the full symmetric
group, or alternating} \label{small} \lremind{small}

\epoint{Derksen's example in $G(4,8)$} 
One of the 31 examples in $G(4,8)$ described above has a Galois group that
is {\em not} at least alternating (and hence is not two-transitive by
our earlier discussion):  the Schubert problem of Figure~\ref{A}.
This example (and the existence of Schubert problems with non-full
Galois group) is due to H.~Derksen.
By Theorem~\ref{galeg}, this is the smallest example of a Schubert
problem with a Galois group smaller than alternating.

\begin{figure}[ht]
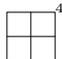

\begin{center}
\include{A}
\end{center}
\caption{This Schubert problem (in $G(4,8)$) has $6$ solutions; the Galois group is $S_4$
\lremind{A}}
\label{A}
\end{figure}

The problem has six solutions.  We show now that the Galois group is
$S_4$.  Fix four general flags in $K^8$, and consider the Schubert
problem in $G(2,8)$ (corresponding to these flags) shown in
Figure~\ref{B}.  This problem has four solutions, corresponding to
four transverse $2$-planes $V_1$, $V_2$, $V_3$, $V_4$.  It is
straightforward to check that $V_i+V_j$ ($i<j$) is a solution to the
original problem of Figure~\ref{A}, in $G(4,8)$.  Hence the Galois
group of the original problem is not two-transitive: two solutions
$W_1$ and $W_2$ may have intersection of dimension $0$ or $2$, and
both possibilities occur.  The Galois group is a subgroup of $S_4$
(acting on the six elements as described above), and is canonically
the Galois group of the problem of Figure~\ref{B}.

\begin{figure}[ht]
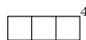

\begin{center}
\include{B}
\end{center}
\caption{An auxiliary Schubert problem in $G(2,8)$ \lremind{B}}
\label{B}
\end{figure}

Applying Theorem~\ref{monthm} (indeed Theorem~\ref{galeg}), the Galois
group of Figure~\ref{B} is at least $A_4$.  By examining the directed
tree of Theorem~\ref{monthm} more closely, we see that the Galois group
is actually $S_4$:  the first branching has one branch with three
leaves and one branch with one leaf (see Remark~\ref{smallrk}).

\epoint{A family of examples generalizing Derksen's} \label{family1}
\lremind{family1} Derksen's example can be generalized to produce
other examples of smaller-than-expected Galois groups, where the
Galois action is that of $S_N$ acting on the order $K$ subsets of $\{
1, 2, \dots, N \}$, as follows.  The Schubert problem of
Figure~\ref{C} in $G(2K, 2N)$ has $\binom N K$ solutions.  Given four
general flags in $G(2,2N)$, the auxiliary problem of Figure~\ref{D}
has $N$ solutions, corresponding to $N$ transverse $2$-planes $V_1$,
\dots, $V_N$ in $G(2,2N)$. By repeated applications of
Remark~\ref{smallrk}, the Galois group of the auxiliary Schubert
problem is $S_N$.  The subspace $V_{i_1} + \cdots + V_{i_K}$ ($1 \leq
i_1 < \cdots < i_K \leq N$) is a solution to the original problem of
Figure~\ref{C}.  Hence the original problem exhibits the desired
behavior.

\begin{figure}[ht]
\begin{center}
\include{C}
\end{center}
\caption{A Schubert problem in $G(2K,2N)$ with $\binom N K$ solutions
and Galois group $S_N$ \lremind{C}}
\label{C}
\end{figure}

\begin{figure}[ht]
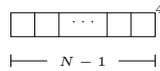

\begin{center}
\include{D}
\end{center}
\caption{An auxiliary problem \lremind{D}}
\label{D}
\end{figure}

The only statements in the previous paragraph that are nontrivial to
verify are (i) the enumeration of solutions to the Schubert problem,
and (ii) the fact that the Galois group of the auxiliary problem is
$S_N$.  Both are easiest to see in terms of puzzles.  (See \cite{ktw}
for a definition of puzzles, and the appendix to \cite{geolr} for the
bijection between checkers and puzzles.)

Part (i) is the number of ways of filling in the puzzle of
Figure~\ref{E} (where the blocks of $1$'s are all of size $K$, and the
blocks of $0$'s are all of size $N-K$), which reduces to
Figure~\ref{F}.  After trying the puzzle, the reader will quickly see
that the number of solutions is $\binom N K$.  The solutions
correspond to the choice of labels on segment $A$ --- there will be
$N-K$ 0's and $K$ 1's, and each order appears in precisely one
completed puzzle.

\begin{figure}[ht]
\begin{center}
\include{E}
\end{center}
\caption{The puzzle corresponding to Figure~\ref{C} \lremind{E}}
\label{E}
\end{figure}

\begin{figure}[ht]
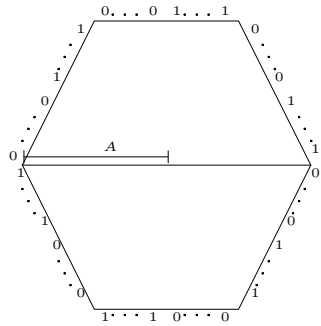

\begin{center}
\include{F}
\end{center}
\caption{The puzzle corresponding to Figure~\ref{C}, partially completed
\lremind{F}}
\label{F}
\end{figure}

To construct the directed tree for part (ii), note that the order of the first
checker game corresponds to filling in the top half of the puzzle of
Figure~\ref{E} (and hence Figure~\ref{F}) row by row; the directed
graph corresponds to the tree of choices made while completing the
puzzle in this order.
Applying this in the case $K=1$, the puzzles of the previous paragraph
show that the tree is of the desired form.

It is interesting (but inessential) to note more generally that the
tree for $\binom N K$ (call it of type $(N,K)$) can be interpreted in
terms of Pascal's triangle as follows.  The two branches at the first
branch point have $\binom {N-1} {K-1}$ and $\binom {N-1} K$ leaves,
and the two corresponding directed trees are of type $(N-1, K-1)$ and
$(N-1, K)$ respectively.  Thus Theorem~\ref{monthm} fails to apply
because of vertices of type $(2 N'', N'')$, corresponding to the
central terms in Pascal's triangle.

\epoint{A similar family of three-flag examples} \label{family2}
\lremind{family2}
We now exhibit a family of three-flag examples with behavior similar to that
of the previous section.    The Schubert problem of Figure~\ref{G}
in $G(3K, 3N)$
has $\binom N K$ solutions and Galois group $S_N$, where the action
is that of $S_N$ on order $K$ subsets of $\{ 1, \dots, N \}$.

\begin{figure}[ht]
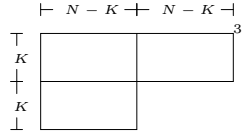

\begin{center}
\include{G}
\end{center}
\caption{A Schubert problem in $G(3K, 3N)$ with $\binom N K$ solutions
and Galois group $S_N$
\lremind{G}}
\label{G}
\end{figure}

As with the previous family, to prove this, first count solutions
using checkers or puzzles.  The puzzle is shown in Figure~\ref{H},
which again reduces to Figure~\ref{F} (without the equatorial cut).
Next, fix three general flags.  Consider the analogous problem with
$K=1$.  There are $N$ solutions, corresponding to $N$ transverse
3-spaces $V_1$, \dots, $V_N$.  The Galois group is $S_N$ by
Remark~\ref{smallrk}, as the tree is identical to that of the previous
section.  The sum of any $K$ of these 3-spaces is a solution to the
original Schubert problem (with respect to the same three flags).
Thus the Galois group of the original problem is $S_N$ as desired.

\begin{figure}[ht]
\begin{center}
\include{H}
\end{center}
\caption{The puzzle corresponding to Figure~\ref{G}
\lremind{H}}
\label{H}
\end{figure}

} 

\end{document}

%% file: g36.tex
\setlength{\unitlength}{0.00083333in}
\begingroup\makeatletter\ifx\SetFigFont\undefined%
\gdef\SetFigFont#1#2#3#4#5{%
  \reset@font\fontsize{#1}{#2pt}%
  \fontfamily{#3}\fontseries{#4}\fontshape{#5}%
  \selectfont}%
\fi\endgroup%
{\renewcommand{\dashlinestretch}{30}
\begin{picture}(1299,411)(0,-10)
\path(12,312)(162,312)(162,162)
	(12,162)(12,312)
\path(612,312)(912,312)(912,12)
	(612,12)(612,312)
\path(1137,312)(1287,312)(1287,162)
	(1137,162)(1137,312)
\path(762,312)(762,12)
\path(612,162)(912,162)
\put(162,312){\makebox(0,0)[lb]{\smash{{{\SetFigFont{5}{6.0}{\rmdefault}{\mddefault}{\updefault}$9$}}}}}
\put(1287,312){\makebox(0,0)[lb]{\smash{{{\SetFigFont{5}{6.0}{\rmdefault}{\mddefault}{\updefault}$5$}}}}}
\put(987,237){\makebox(0,0)[lb]{\smash{{{\SetFigFont{5}{6.0}{\rmdefault}{\mddefault}{\updefault}$\times$}}}}}
\end{picture}
}

%% file: g37a.tex
\setlength{\unitlength}{0.00083333in}
\begingroup\makeatletter\ifx\SetFigFont\undefined%
\gdef\SetFigFont#1#2#3#4#5{%
  \reset@font\fontsize{#1}{#2pt}%
  \fontfamily{#3}\fontseries{#4}\fontshape{#5}%
  \selectfont}%
\fi\endgroup%
{\renewcommand{\dashlinestretch}{30}
\begin{picture}(1224,561)(0,-10)
\path(12,462)(312,462)(312,162)
	(12,162)(12,462)
\path(537,462)(687,462)(687,312)
	(537,312)(537,462)
\path(912,462)(1062,462)(1062,12)
	(912,12)(912,462)
\path(912,312)(1212,312)(1212,462)
	(912,462)(912,312)
\path(912,162)(1062,162)
\path(12,312)(312,312)
\path(162,462)(162,162)
\put(387,387){\makebox(0,0)[lb]{\smash{{{\SetFigFont{5}{6.0}{\rmdefault}{\mddefault}{\updefault}$\times$}}}}}
\put(762,387){\makebox(0,0)[lb]{\smash{{{\SetFigFont{5}{6.0}{\rmdefault}{\mddefault}{\updefault}$\times$}}}}}
\put(687,462){\makebox(0,0)[lb]{\smash{{{\SetFigFont{5}{6.0}{\rmdefault}{\mddefault}{\updefault}$4$}}}}}
\end{picture}
}

%% file: g37.tex
\setlength{\unitlength}{0.00083333in}
\begingroup\makeatletter\ifx\SetFigFont\undefined%
\gdef\SetFigFont#1#2#3#4#5{%
  \reset@font\fontsize{#1}{#2pt}%
  \fontfamily{#3}\fontseries{#4}\fontshape{#5}%
  \selectfont}%
\fi\endgroup%
{\renewcommand{\dashlinestretch}{30}
\begin{picture}(3774,411)(0,-10)
\path(2187,312)(2637,312)(2637,12)
	(2187,12)(2187,312)
\path(2337,312)(2337,12)
\path(2487,312)(2487,12)
\path(2187,162)(2637,162)
\path(1287,312)(1437,312)(1437,162)
	(1287,162)(1287,312)
\path(1662,312)(1962,312)(1962,162)
	(1662,162)(1662,312)
\path(1812,312)(1812,162)
\put(1437,312){\makebox(0,0)[lb]{\smash{{{\SetFigFont{5}{6.0}{\rmdefault}{\mddefault}{\updefault}$4$}}}}}
\put(1512,237){\makebox(0,0)[lb]{\smash{{{\SetFigFont{5}{6.0}{\rmdefault}{\mddefault}{\updefault}$\times$}}}}}
\put(2037,237){\makebox(0,0)[lb]{\smash{{{\SetFigFont{5}{6.0}{\rmdefault}{\mddefault}{\updefault}$\times$}}}}}
\path(387,312)(837,312)(837,12)
	(387,12)(387,312)
\path(537,312)(537,12)
\path(687,312)(687,12)
\path(387,162)(837,162)
\path(12,312)(162,312)(162,162)
	(12,162)(12,312)
\put(162,312){\makebox(0,0)[lb]{\smash{{{\SetFigFont{5}{6.0}{\rmdefault}{\mddefault}{\updefault}$6$}}}}}
\put(237,237){\makebox(0,0)[lb]{\smash{{{\SetFigFont{5}{6.0}{\rmdefault}{\mddefault}{\updefault}$\times$}}}}}
\path(3087,12)(3237,12)(3237,312)
	(3087,312)(3087,12)
\path(3462,312)(3762,312)(3762,162)
	(3462,162)(3462,312)
\path(3087,162)(3237,162)
\path(3612,312)(3612,162)
\put(3312,237){\makebox(0,0)[lb]{\smash{{{\SetFigFont{5}{6.0}{\rmdefault}{\mddefault}{\updefault}$\times$}}}}}
\put(3237,312){\makebox(0,0)[lb]{\smash{{{\SetFigFont{5}{6.0}{\rmdefault}{\mddefault}{\updefault}$5$}}}}}
\end{picture}
}

%% file: g38.tex
\setlength{\unitlength}{0.00083333in}
\begingroup\makeatletter\ifx\SetFigFont\undefined%
\gdef\SetFigFont#1#2#3#4#5{%
  \reset@font\fontsize{#1}{#2pt}%
  \fontfamily{#3}\fontseries{#4}\fontshape{#5}%
  \selectfont}%
\fi\endgroup%
{\renewcommand{\dashlinestretch}{30}
\begin{picture}(5574,1161)(0,-10)
\path(387,1062)(837,1062)(837,912)
	(387,912)(387,1062)
\path(537,1062)(537,912)
\path(687,1062)(687,912)
\path(1062,1062)(1662,1062)(1662,762)
	(1062,762)(1062,1062)
\path(1212,1062)(1212,762)
\path(1362,1062)(1362,762)
\path(1512,1062)(1512,762)
\path(1062,912)(1662,912)
\path(2487,1062)(2787,1062)(2787,912)
	(2487,912)(2487,1062)
\path(2637,1062)(2637,912)
\path(3012,1062)(3612,1062)(3612,762)
	(3012,762)(3012,1062)
\path(3162,1062)(3162,762)
\path(3312,1062)(3312,762)
\path(3462,1062)(3462,762)
\path(3012,912)(3612,912)
\path(4437,1062)(4737,1062)(4737,912)
	(4437,912)(4437,1062)
\path(4587,1062)(4587,912)
\path(4962,1062)(5562,1062)(5562,762)
	(4962,762)(4962,1062)
\path(5112,1062)(5112,762)
\path(5262,1062)(5262,762)
\path(5412,1062)(5412,762)
\path(4962,912)(5562,912)
\path(387,312)(987,312)(987,12)
	(387,12)(387,312)
\path(537,312)(537,12)
\path(687,312)(687,12)
\path(837,312)(837,12)
\path(387,162)(987,162)
\path(2112,312)(2262,312)(2262,12)
	(2112,12)(2112,312)
\path(2112,162)(2262,162)
\path(4062,312)(4212,312)(4212,12)
	(4062,12)(4062,312)
\path(4062,162)(4212,162)
\path(4437,312)(4887,312)(4887,162)
	(4437,162)(4437,312)
\path(4587,312)(4587,162)
\path(4737,312)(4737,162)
\path(12,1062)(162,1062)(162,912)
	(12,912)(12,1062)
\path(2112,1062)(2262,1062)(2262,912)
	(2112,912)(2112,1062)
\path(4062,1062)(4212,1062)(4212,912)
	(4062,912)(4062,1062)
\path(12,312)(162,312)(162,162)
	(12,162)(12,312)
\path(2487,312)(2787,312)(2787,12)
	(2487,12)(2487,312)
\path(3012,312)(3462,312)(3462,162)
	(3012,162)(3012,312)
\path(3162,312)(3162,162)
\path(3312,312)(3312,162)
\path(2487,162)(2787,162)
\path(2637,312)(2637,12)
\put(162,1062){\makebox(0,0)[lb]{\smash{{{\SetFigFont{5}{6.0}{\rmdefault}{\mddefault}{\updefault}$4$}}}}}
\put(237,987){\makebox(0,0)[lb]{\smash{{{\SetFigFont{5}{6.0}{\rmdefault}{\mddefault}{\updefault}$\times$}}}}}
\put(912,987){\makebox(0,0)[lb]{\smash{{{\SetFigFont{5}{6.0}{\rmdefault}{\mddefault}{\updefault}$\times$}}}}}
\put(2337,987){\makebox(0,0)[lb]{\smash{{{\SetFigFont{5}{6.0}{\rmdefault}{\mddefault}{\updefault}$\times$}}}}}
\put(2862,987){\makebox(0,0)[lb]{\smash{{{\SetFigFont{5}{6.0}{\rmdefault}{\mddefault}{\updefault}$\times$}}}}}
\put(2262,1062){\makebox(0,0)[lb]{\smash{{{\SetFigFont{5}{6.0}{\rmdefault}{\mddefault}{\updefault}$3$}}}}}
\put(2787,1062){\makebox(0,0)[lb]{\smash{{{\SetFigFont{5}{6.0}{\rmdefault}{\mddefault}{\updefault}$2$}}}}}
\put(4287,987){\makebox(0,0)[lb]{\smash{{{\SetFigFont{5}{6.0}{\rmdefault}{\mddefault}{\updefault}$\times$}}}}}
\put(4812,987){\makebox(0,0)[lb]{\smash{{{\SetFigFont{5}{6.0}{\rmdefault}{\mddefault}{\updefault}$\times$}}}}}
\put(4212,1062){\makebox(0,0)[lb]{\smash{{{\SetFigFont{5}{6.0}{\rmdefault}{\mddefault}{\updefault}$5$}}}}}
\put(237,237){\makebox(0,0)[lb]{\smash{{{\SetFigFont{5}{6.0}{\rmdefault}{\mddefault}{\updefault}$\times$}}}}}
\put(162,312){\makebox(0,0)[lb]{\smash{{{\SetFigFont{5}{6.0}{\rmdefault}{\mddefault}{\updefault}$7$}}}}}
\put(2262,312){\makebox(0,0)[lb]{\smash{{{\SetFigFont{5}{6.0}{\rmdefault}{\mddefault}{\updefault}$4$}}}}}
\put(2337,237){\makebox(0,0)[lb]{\smash{{{\SetFigFont{5}{6.0}{\rmdefault}{\mddefault}{\updefault}$\times$}}}}}
\put(2862,237){\makebox(0,0)[lb]{\smash{{{\SetFigFont{5}{6.0}{\rmdefault}{\mddefault}{\updefault}$\times$}}}}}
\put(4212,312){\makebox(0,0)[lb]{\smash{{{\SetFigFont{5}{6.0}{\rmdefault}{\mddefault}{\updefault}$6$}}}}}
\put(4287,237){\makebox(0,0)[lb]{\smash{{{\SetFigFont{5}{6.0}{\rmdefault}{\mddefault}{\updefault}$\times$}}}}}
\end{picture}
}

%% file: g39.tex
\setlength{\unitlength}{0.00083333in}
\begingroup\makeatletter\ifx\SetFigFont\undefined%
\gdef\SetFigFont#1#2#3#4#5{%
  \reset@font\fontsize{#1}{#2pt}%
  \fontfamily{#3}\fontseries{#4}\fontshape{#5}%
  \selectfont}%
\fi\endgroup%
{\renewcommand{\dashlinestretch}{30}
\begin{picture}(6024,3411)(0,-10)
\path(387,3312)(987,3312)(987,3162)
	(387,3162)(387,3312)
\path(537,3312)(537,3162)
\path(687,3312)(687,3162)
\path(837,3312)(837,3162)
\path(1212,3312)(1962,3312)(1962,3012)
	(1212,3012)(1212,3312)
\path(1212,3162)(1962,3162)
\path(1512,3312)(1512,3012)
\path(1662,3312)(1662,3012)
\path(1812,3312)(1812,3012)
\path(1362,3312)(1362,3012)
\path(2787,3312)(3087,3312)(3087,3162)
	(2787,3162)(2787,3312)
\path(2937,3312)(2937,3162)
\path(3312,3312)(3762,3312)(3762,3162)
	(3312,3162)(3312,3312)
\path(3462,3312)(3462,3162)
\path(3612,3312)(3612,3162)
\path(3987,3312)(4737,3312)(4737,3012)
	(3987,3012)(3987,3312)
\path(3987,3162)(4737,3162)
\path(4287,3312)(4287,3012)
\path(4437,3312)(4437,3012)
\path(4587,3312)(4587,3012)
\path(4137,3312)(4137,3012)
\path(387,2562)(837,2562)(837,2412)
	(387,2412)(387,2562)
\path(537,2562)(537,2412)
\path(687,2562)(687,2412)
\path(1062,2562)(1812,2562)(1812,2262)
	(1062,2262)(1062,2562)
\path(1062,2412)(1812,2412)
\path(1362,2562)(1362,2262)
\path(1512,2562)(1512,2262)
\path(1662,2562)(1662,2262)
\path(1212,2562)(1212,2262)
\path(2637,2562)(2937,2562)(2937,2412)
	(2637,2412)(2637,2562)
\path(2787,2562)(2787,2412)
\path(3162,2562)(3912,2562)(3912,2262)
	(3162,2262)(3162,2562)
\path(3162,2412)(3912,2412)
\path(3462,2562)(3462,2262)
\path(3612,2562)(3612,2262)
\path(3762,2562)(3762,2262)
\path(3312,2562)(3312,2262)
\path(387,1812)(687,1812)(687,1662)
	(387,1662)(387,1812)
\path(537,1812)(537,1662)
\path(912,1812)(1662,1812)(1662,1512)
	(912,1512)(912,1812)
\path(912,1662)(1662,1662)
\path(1212,1812)(1212,1512)
\path(1362,1812)(1362,1512)
\path(1512,1812)(1512,1512)
\path(1062,1812)(1062,1512)
\path(4737,2562)(5037,2562)(5037,2412)
	(4737,2412)(4737,2562)
\path(4887,2562)(4887,2412)
\path(5262,2562)(6012,2562)(6012,2262)
	(5262,2262)(5262,2562)
\path(5262,2412)(6012,2412)
\path(5562,2562)(5562,2262)
\path(5712,2562)(5712,2262)
\path(5862,2562)(5862,2262)
\path(5412,2562)(5412,2262)
\path(2637,1812)(3387,1812)(3387,1512)
	(2637,1512)(2637,1812)
\path(2637,1662)(3387,1662)
\path(2937,1812)(2937,1512)
\path(3087,1812)(3087,1512)
\path(3237,1812)(3237,1512)
\path(2787,1812)(2787,1512)
\path(4362,1812)(4962,1812)(4962,1662)
	(4362,1662)(4362,1812)
\path(4512,1812)(4512,1662)
\path(4662,1812)(4662,1662)
\path(4812,1812)(4812,1662)
\path(5187,1812)(5337,1812)(5337,1512)
	(5187,1512)(5187,1812)
\path(5187,1662)(5337,1662)
\path(5562,1812)(6012,1812)(6012,1512)
	(5562,1512)(5562,1812)
\path(5562,1662)(6012,1662)
\path(5862,1812)(5862,1512)
\path(5712,1812)(5712,1512)
\path(12,1062)(462,1062)(462,912)
	(12,912)(12,1062)
\path(162,1062)(162,912)
\path(312,1062)(312,912)
\path(687,1062)(1287,1062)(1287,912)
	(687,912)(687,1062)
\path(837,1062)(837,912)
\path(987,1062)(987,912)
\path(1137,1062)(1137,912)
\path(2412,1062)(2862,1062)(2862,912)
	(2412,912)(2412,1062)
\path(2562,1062)(2562,912)
\path(2712,1062)(2712,912)
\path(3087,1062)(3237,1062)(3237,762)
	(3087,762)(3087,1062)
\path(3087,912)(3237,912)
\path(3462,1062)(3762,1062)(3762,762)
	(3462,762)(3462,1062)
\path(3462,912)(3762,912)
\path(3612,1062)(3612,762)
\put(2937,987){\makebox(0,0)[lb]{\smash{{{\SetFigFont{5}{6.0}{\rmdefault}{\mddefault}{\updefault}$\times$}}}}}
\put(3312,987){\makebox(0,0)[lb]{\smash{{{\SetFigFont{5}{6.0}{\rmdefault}{\mddefault}{\updefault}$\times$}}}}}
\put(3237,1062){\makebox(0,0)[lb]{\smash{{{\SetFigFont{5}{6.0}{\rmdefault}{\mddefault}{\updefault}$3$}}}}}
\put(3762,1062){\makebox(0,0)[lb]{\smash{{{\SetFigFont{5}{6.0}{\rmdefault}{\mddefault}{\updefault}$2$}}}}}
\path(2412,1062)(2262,1062)(2262,912)(2412,912)
\path(4512,1062)(4962,1062)(4962,912)
	(4512,912)(4512,1062)
\path(4662,1062)(4662,912)
\path(4812,1062)(4812,912)
\path(5187,1062)(5337,1062)(5337,762)
	(5187,762)(5187,1062)
\path(5187,912)(5337,912)
\path(5562,1062)(5862,1062)(5862,762)
	(5562,762)(5562,1062)
\path(5562,912)(5862,912)
\path(5712,1062)(5712,762)
\path(4512,1062)(4362,1062)(4362,912)(4512,912)
\put(5037,987){\makebox(0,0)[lb]{\smash{{{\SetFigFont{5}{6.0}{\rmdefault}{\mddefault}{\updefault}$\times$}}}}}
\put(5412,987){\makebox(0,0)[lb]{\smash{{{\SetFigFont{5}{6.0}{\rmdefault}{\mddefault}{\updefault}$\times$}}}}}
\put(5337,1062){\makebox(0,0)[lb]{\smash{{{\SetFigFont{5}{6.0}{\rmdefault}{\mddefault}{\updefault}$5$}}}}}
\path(162,312)(612,312)(612,162)
	(162,162)(162,312)
\path(312,312)(312,162)
\path(462,312)(462,162)
\path(2262,312)(2412,312)(2412,12)
	(2262,12)(2262,312)
\path(2262,162)(2412,162)
\path(837,312)(987,312)(987,12)
	(837,12)(837,312)
\path(837,162)(987,162)
\path(12,3312)(162,3312)(162,3162)
	(12,3162)(12,3312)
\path(2412,3312)(2562,3312)(2562,3162)
	(2412,3162)(2412,3312)
\path(12,2562)(162,2562)(162,2412)
	(12,2412)(12,2562)
\path(2262,2562)(2412,2562)(2412,2412)
	(2262,2412)(2262,2562)
\path(4362,2562)(4512,2562)(4512,2412)
	(4362,2412)(4362,2562)
\path(12,1812)(162,1812)(162,1662)
	(12,1662)(12,1812)
\path(2262,1812)(2412,1812)(2412,1662)
	(2262,1662)(2262,1812)
\path(687,912)(687,762)(837,762)
	(837,912)(687,912)
\path(162,312)(12,312)(12,162)(162,162)
\put(237,3237){\makebox(0,0)[lb]{\smash{{{\SetFigFont{5}{6.0}{\rmdefault}{\mddefault}{\updefault}$\times$}}}}}
\put(1062,3237){\makebox(0,0)[lb]{\smash{{{\SetFigFont{5}{6.0}{\rmdefault}{\mddefault}{\updefault}$\times$}}}}}
\put(162,3312){\makebox(0,0)[lb]{\smash{{{\SetFigFont{5}{6.0}{\rmdefault}{\mddefault}{\updefault}$4$}}}}}
\put(2637,3237){\makebox(0,0)[lb]{\smash{{{\SetFigFont{5}{6.0}{\rmdefault}{\mddefault}{\updefault}$\times$}}}}}
\put(3162,3237){\makebox(0,0)[lb]{\smash{{{\SetFigFont{5}{6.0}{\rmdefault}{\mddefault}{\updefault}$\times$}}}}}
\put(3837,3237){\makebox(0,0)[lb]{\smash{{{\SetFigFont{5}{6.0}{\rmdefault}{\mddefault}{\updefault}$\times$}}}}}
\put(2562,3312){\makebox(0,0)[lb]{\smash{{{\SetFigFont{5}{6.0}{\rmdefault}{\mddefault}{\updefault}$3$}}}}}
\put(162,2562){\makebox(0,0)[lb]{\smash{{{\SetFigFont{5}{6.0}{\rmdefault}{\mddefault}{\updefault}$5$}}}}}
\put(237,2487){\makebox(0,0)[lb]{\smash{{{\SetFigFont{5}{6.0}{\rmdefault}{\mddefault}{\updefault}$\times$}}}}}
\put(912,2487){\makebox(0,0)[lb]{\smash{{{\SetFigFont{5}{6.0}{\rmdefault}{\mddefault}{\updefault}$\times$}}}}}
\put(2487,2487){\makebox(0,0)[lb]{\smash{{{\SetFigFont{5}{6.0}{\rmdefault}{\mddefault}{\updefault}$\times$}}}}}
\put(3012,2487){\makebox(0,0)[lb]{\smash{{{\SetFigFont{5}{6.0}{\rmdefault}{\mddefault}{\updefault}$\times$}}}}}
\put(2937,2562){\makebox(0,0)[lb]{\smash{{{\SetFigFont{5}{6.0}{\rmdefault}{\mddefault}{\updefault}$3$}}}}}
\put(2412,2562){\makebox(0,0)[lb]{\smash{{{\SetFigFont{5}{6.0}{\rmdefault}{\mddefault}{\updefault}$2$}}}}}
\put(237,1737){\makebox(0,0)[lb]{\smash{{{\SetFigFont{5}{6.0}{\rmdefault}{\mddefault}{\updefault}$\times$}}}}}
\put(762,1737){\makebox(0,0)[lb]{\smash{{{\SetFigFont{5}{6.0}{\rmdefault}{\mddefault}{\updefault}$\times$}}}}}
\put(4587,2487){\makebox(0,0)[lb]{\smash{{{\SetFigFont{5}{6.0}{\rmdefault}{\mddefault}{\updefault}$\times$}}}}}
\put(5112,2487){\makebox(0,0)[lb]{\smash{{{\SetFigFont{5}{6.0}{\rmdefault}{\mddefault}{\updefault}$\times$}}}}}
\put(5037,2562){\makebox(0,0)[lb]{\smash{{{\SetFigFont{5}{6.0}{\rmdefault}{\mddefault}{\updefault}$2$}}}}}
\put(4512,2562){\makebox(0,0)[lb]{\smash{{{\SetFigFont{5}{6.0}{\rmdefault}{\mddefault}{\updefault}$4$}}}}}
\put(162,1812){\makebox(0,0)[lb]{\smash{{{\SetFigFont{5}{6.0}{\rmdefault}{\mddefault}{\updefault}$6$}}}}}
\put(2487,1737){\makebox(0,0)[lb]{\smash{{{\SetFigFont{5}{6.0}{\rmdefault}{\mddefault}{\updefault}$\times$}}}}}
\put(2412,1812){\makebox(0,0)[lb]{\smash{{{\SetFigFont{5}{6.0}{\rmdefault}{\mddefault}{\updefault}$8$}}}}}
\put(5412,1737){\makebox(0,0)[lb]{\smash{{{\SetFigFont{5}{6.0}{\rmdefault}{\mddefault}{\updefault}$\times$}}}}}
\put(5037,1737){\makebox(0,0)[lb]{\smash{{{\SetFigFont{5}{6.0}{\rmdefault}{\mddefault}{\updefault}$\times$}}}}}
\put(5337,1812){\makebox(0,0)[lb]{\smash{{{\SetFigFont{5}{6.0}{\rmdefault}{\mddefault}{\updefault}$4$}}}}}
\put(537,987){\makebox(0,0)[lb]{\smash{{{\SetFigFont{5}{6.0}{\rmdefault}{\mddefault}{\updefault}$\times$}}}}}
\put(1287,1062){\makebox(0,0)[lb]{\smash{{{\SetFigFont{5}{6.0}{\rmdefault}{\mddefault}{\updefault}$3$}}}}}
\put(687,237){\makebox(0,0)[lb]{\smash{{{\SetFigFont{5}{6.0}{\rmdefault}{\mddefault}{\updefault}$\times$}}}}}
\put(987,312){\makebox(0,0)[lb]{\smash{{{\SetFigFont{5}{6.0}{\rmdefault}{\mddefault}{\updefault}$7$}}}}}
\put(2412,312){\makebox(0,0)[lb]{\smash{{{\SetFigFont{5}{6.0}{\rmdefault}{\mddefault}{\updefault}$9$}}}}}
\end{picture}
}

%% file: A.tex
\setlength{\unitlength}{0.00083333in}
\begingroup\makeatletter\ifx\SetFigFont\undefined%
\gdef\SetFigFont#1#2#3#4#5{%
  \reset@font\fontsize{#1}{#2pt}%
  \fontfamily{#3}\fontseries{#4}\fontshape{#5}%
  \selectfont}%
\fi\endgroup%
{\renewcommand{\dashlinestretch}{30}
\begin{picture}(324,411)(0,-10)
\path(12,312)(312,312)(312,12)
	(12,12)(12,312)
\path(162,312)(162,12)
\path(12,162)(312,162)
\put(312,312){\makebox(0,0)[lb]{\smash{{{\SetFigFont{5}{6.0}{\rmdefault}{\mddefault}{\updefault}$4$}}}}}
\end{picture}
}

%% file: B.tex
\setlength{\unitlength}{0.00083333in}
\begingroup\makeatletter\ifx\SetFigFont\undefined%
\gdef\SetFigFont#1#2#3#4#5{%
  \reset@font\fontsize{#1}{#2pt}%
  \fontfamily{#3}\fontseries{#4}\fontshape{#5}%
  \selectfont}%
\fi\endgroup%
{\renewcommand{\dashlinestretch}{30}
\begin{picture}(474,261)(0,-10)
\path(462,162)(12,162)(12,12)
	(462,12)(462,162)
\path(162,162)(162,12)
\path(312,162)(312,12)
\put(462,162){\makebox(0,0)[lb]{\smash{{{\SetFigFont{5}{6.0}{\rmdefault}{\mddefault}{\updefault}$4$}}}}}
\end{picture}
}

%% file: C.tex
\setlength{\unitlength}{0.00083333in}
\begingroup\makeatletter\ifx\SetFigFont\undefined%
\gdef\SetFigFont#1#2#3#4#5{%
  \reset@font\fontsize{#1}{#2pt}%
  \fontfamily{#3}\fontseries{#4}\fontshape{#5}%
  \selectfont}%
\fi\endgroup%
{\renewcommand{\dashlinestretch}{30}
\begin{picture}(1112,910)(0,-10)
\path(875,61)(1100,61)
\path(200,61)(425,61)
\path(200,99)(200,24)
\path(1100,99)(1100,24)
\path(1100,811)(200,811)(200,211)
	(1100,211)(1100,811)
\path(500,811)(500,661)(200,661)
\path(350,811)(350,661)
\path(950,811)(950,661)
\path(800,811)(800,661)(1100,661)
\path(200,361)(500,361)(500,211)
\path(800,211)(800,361)(1100,361)
\path(350,361)(350,211)
\path(950,361)(950,211)
\path(50,811)(50,586)
\path(50,436)(50,211)
\path(12,811)(87,811)
\path(12,211)(87,211)
\put(500,31){\makebox(0,0)[lb]{\smash{{{\SetFigFont{5}{6.0}{\rmdefault}{\mddefault}{\updefault}$N-K$}}}}}
\put(1100,811){\makebox(0,0)[lb]{\smash{{{\SetFigFont{5}{6.0}{\rmdefault}{\mddefault}{\updefault}$4$}}}}}
\put(30,474){\makebox(0,0)[lb]{\smash{{{\SetFigFont{5}{6.0}{\rmdefault}{\mddefault}{\updefault}$K$}}}}}
\put(575,474){\makebox(0,0)[lb]{\smash{{{\SetFigFont{5}{6.0}{\rmdefault}{\mddefault}{\updefault}$\ddots$}}}}}
\end{picture}
}

%% file: D.tex
\setlength{\unitlength}{0.00083333in}
\begingroup\makeatletter\ifx\SetFigFont\undefined%
\gdef\SetFigFont#1#2#3#4#5{%
  \reset@font\fontsize{#1}{#2pt}%
  \fontfamily{#3}\fontseries{#4}\fontshape{#5}%
  \selectfont}%
\fi\endgroup%
{\renewcommand{\dashlinestretch}{30}
\begin{picture}(924,460)(0,-10)
\path(612,361)(612,211)
\path(762,361)(762,211)
\path(912,361)(12,361)(12,211)
	(912,211)(912,361)
\path(162,361)(162,211)
\path(312,361)(312,211)
\path(687,61)(912,61)
\path(12,61)(237,61)
\path(12,99)(12,24)
\path(912,99)(912,24)
\put(912,361){\makebox(0,0)[lb]{\smash{{{\SetFigFont{5}{6.0}{\rmdefault}{\mddefault}{\updefault}$4$}}}}}
\put(387,286){\makebox(0,0)[lb]{\smash{{{\SetFigFont{5}{6.0}{\rmdefault}{\mddefault}{\updefault}$\cdots$}}}}}
\put(312,31){\makebox(0,0)[lb]{\smash{{{\SetFigFont{5}{6.0}{\rmdefault}{\mddefault}{\updefault}$N-1$}}}}}
\end{picture}
}

%% file: E.tex
\setlength{\unitlength}{0.00083333in}
\begingroup\makeatletter\ifx\SetFigFont\undefined%
\gdef\SetFigFont#1#2#3#4#5{%
  \reset@font\fontsize{#1}{#2pt}%
  \fontfamily{#3}\fontseries{#4}\fontshape{#5}%
  \selectfont}%
\fi\endgroup%
{\renewcommand{\dashlinestretch}{30}
\begin{picture}(5325,3638)(0,-10)
\put(2662,2561){\blacken\ellipse{10}{10}}
\put(2662,2561){\ellipse{10}{10}}
\put(2700,2486){\blacken\ellipse{10}{10}}
\put(2700,2486){\ellipse{10}{10}}
\put(2737,2411){\blacken\ellipse{10}{10}}
\put(2737,2411){\ellipse{10}{10}}
\put(2925,2036){\blacken\ellipse{10}{10}}
\put(2925,2036){\ellipse{10}{10}}
\put(2963,1961){\blacken\ellipse{10}{10}}
\put(2963,1961){\ellipse{10}{10}}
\put(2888,2111){\blacken\ellipse{10}{10}}
\put(2888,2111){\ellipse{10}{10}}
\put(2587,2636){\makebox(0,0)[lb]{\smash{{{\SetFigFont{5}{6.0}{\rmdefault}{\mddefault}{\updefault}$0$}}}}}
\put(2963,1886){\makebox(0,0)[lb]{\smash{{{\SetFigFont{5}{6.0}{\rmdefault}{\mddefault}{\updefault}$1$}}}}}
\put(2812,2186){\makebox(0,0)[lb]{\smash{{{\SetFigFont{5}{6.0}{\rmdefault}{\mddefault}{\updefault}$1$}}}}}
\put(2737,2336){\makebox(0,0)[lb]{\smash{{{\SetFigFont{5}{6.0}{\rmdefault}{\mddefault}{\updefault}$0$}}}}}
\put(2212,3461){\blacken\ellipse{10}{10}}
\put(2212,3461){\ellipse{10}{10}}
\put(2250,3386){\blacken\ellipse{10}{10}}
\put(2250,3386){\ellipse{10}{10}}
\put(2287,3311){\blacken\ellipse{10}{10}}
\put(2287,3311){\ellipse{10}{10}}
\put(2475,2936){\blacken\ellipse{10}{10}}
\put(2475,2936){\ellipse{10}{10}}
\put(2512,2861){\blacken\ellipse{10}{10}}
\put(2512,2861){\ellipse{10}{10}}
\put(2437,3011){\blacken\ellipse{10}{10}}
\put(2437,3011){\ellipse{10}{10}}
\put(2137,3536){\makebox(0,0)[lb]{\smash{{{\SetFigFont{5}{6.0}{\rmdefault}{\mddefault}{\updefault}$0$}}}}}
\put(2512,2786){\makebox(0,0)[lb]{\smash{{{\SetFigFont{5}{6.0}{\rmdefault}{\mddefault}{\updefault}$1$}}}}}
\put(2362,3086){\makebox(0,0)[lb]{\smash{{{\SetFigFont{5}{6.0}{\rmdefault}{\mddefault}{\updefault}$1$}}}}}
\put(2287,3236){\makebox(0,0)[lb]{\smash{{{\SetFigFont{5}{6.0}{\rmdefault}{\mddefault}{\updefault}$0$}}}}}
\put(1462,2561){\blacken\ellipse{10}{10}}
\put(1462,2561){\ellipse{10}{10}}
\put(1425,2486){\blacken\ellipse{10}{10}}
\put(1425,2486){\ellipse{10}{10}}
\put(1387,2411){\blacken\ellipse{10}{10}}
\put(1387,2411){\ellipse{10}{10}}
\put(1237,2111){\blacken\ellipse{10}{10}}
\put(1237,2111){\ellipse{10}{10}}
\put(1200,2036){\blacken\ellipse{10}{10}}
\put(1200,2036){\ellipse{10}{10}}
\put(1162,1961){\blacken\ellipse{10}{10}}
\put(1162,1961){\ellipse{10}{10}}
\put(1125,1886){\makebox(0,0)[lb]{\smash{{{\SetFigFont{5}{6.0}{\rmdefault}{\mddefault}{\updefault}$0$}}}}}
\put(1500,2636){\makebox(0,0)[lb]{\smash{{{\SetFigFont{5}{6.0}{\rmdefault}{\mddefault}{\updefault}$1$}}}}}
\put(1350,2336){\makebox(0,0)[lb]{\smash{{{\SetFigFont{5}{6.0}{\rmdefault}{\mddefault}{\updefault}$1$}}}}}
\put(1275,2186){\makebox(0,0)[lb]{\smash{{{\SetFigFont{5}{6.0}{\rmdefault}{\mddefault}{\updefault}$0$}}}}}
\put(1912,3461){\blacken\ellipse{10}{10}}
\put(1912,3461){\ellipse{10}{10}}
\put(1875,3386){\blacken\ellipse{10}{10}}
\put(1875,3386){\ellipse{10}{10}}
\put(1837,3311){\blacken\ellipse{10}{10}}
\put(1837,3311){\ellipse{10}{10}}
\put(1687,3011){\blacken\ellipse{10}{10}}
\put(1687,3011){\ellipse{10}{10}}
\put(1650,2936){\blacken\ellipse{10}{10}}
\put(1650,2936){\ellipse{10}{10}}
\put(1612,2861){\blacken\ellipse{10}{10}}
\put(1612,2861){\ellipse{10}{10}}
\put(1575,2786){\makebox(0,0)[lb]{\smash{{{\SetFigFont{5}{6.0}{\rmdefault}{\mddefault}{\updefault}$0$}}}}}
\put(1950,3536){\makebox(0,0)[lb]{\smash{{{\SetFigFont{5}{6.0}{\rmdefault}{\mddefault}{\updefault}$1$}}}}}
\put(1800,3236){\makebox(0,0)[lb]{\smash{{{\SetFigFont{5}{6.0}{\rmdefault}{\mddefault}{\updefault}$1$}}}}}
\put(1725,3086){\makebox(0,0)[lb]{\smash{{{\SetFigFont{5}{6.0}{\rmdefault}{\mddefault}{\updefault}$0$}}}}}
\put(1162,1662){\blacken\ellipse{10}{10}}
\put(1162,1662){\ellipse{10}{10}}
\put(1200,1587){\blacken\ellipse{10}{10}}
\put(1200,1587){\ellipse{10}{10}}
\put(1237,1512){\blacken\ellipse{10}{10}}
\put(1237,1512){\ellipse{10}{10}}
\put(1387,1212){\blacken\ellipse{10}{10}}
\put(1387,1212){\ellipse{10}{10}}
\put(1425,1137){\blacken\ellipse{10}{10}}
\put(1425,1137){\ellipse{10}{10}}
\put(1462,1062){\blacken\ellipse{10}{10}}
\put(1462,1062){\ellipse{10}{10}}
\put(1125,1737){\makebox(0,0)[lb]{\smash{{{\SetFigFont{5}{6.0}{\rmdefault}{\mddefault}{\updefault}$1$}}}}}
\put(1350,1287){\makebox(0,0)[lb]{\smash{{{\SetFigFont{5}{6.0}{\rmdefault}{\mddefault}{\updefault}$0$}}}}}
\put(1500,987){\makebox(0,0)[lb]{\smash{{{\SetFigFont{5}{6.0}{\rmdefault}{\mddefault}{\updefault}$0$}}}}}
\put(1275,1437){\makebox(0,0)[lb]{\smash{{{\SetFigFont{5}{6.0}{\rmdefault}{\mddefault}{\updefault}$1$}}}}}
\put(1612,762){\blacken\ellipse{10}{10}}
\put(1612,762){\ellipse{10}{10}}
\put(1650,687){\blacken\ellipse{10}{10}}
\put(1650,687){\ellipse{10}{10}}
\put(1687,612){\blacken\ellipse{10}{10}}
\put(1687,612){\ellipse{10}{10}}
\put(1837,312){\blacken\ellipse{10}{10}}
\put(1837,312){\ellipse{10}{10}}
\put(1875,237){\blacken\ellipse{10}{10}}
\put(1875,237){\ellipse{10}{10}}
\put(1912,162){\blacken\ellipse{10}{10}}
\put(1912,162){\ellipse{10}{10}}
\put(1575,837){\makebox(0,0)[lb]{\smash{{{\SetFigFont{5}{6.0}{\rmdefault}{\mddefault}{\updefault}$1$}}}}}
\put(1800,387){\makebox(0,0)[lb]{\smash{{{\SetFigFont{5}{6.0}{\rmdefault}{\mddefault}{\updefault}$0$}}}}}
\put(1950,87){\makebox(0,0)[lb]{\smash{{{\SetFigFont{5}{6.0}{\rmdefault}{\mddefault}{\updefault}$0$}}}}}
\put(1725,537){\makebox(0,0)[lb]{\smash{{{\SetFigFont{5}{6.0}{\rmdefault}{\mddefault}{\updefault}$1$}}}}}
\put(2700,1212){\blacken\ellipse{10}{10}}
\put(2700,1212){\ellipse{10}{10}}
\put(2662,1137){\blacken\ellipse{10}{10}}
\put(2662,1137){\ellipse{10}{10}}
\put(2625,1062){\blacken\ellipse{10}{10}}
\put(2625,1062){\ellipse{10}{10}}
\put(2925,1662){\blacken\ellipse{10}{10}}
\put(2925,1662){\ellipse{10}{10}}
\put(2888,1587){\blacken\ellipse{10}{10}}
\put(2888,1587){\ellipse{10}{10}}
\put(2850,1512){\blacken\ellipse{10}{10}}
\put(2850,1512){\ellipse{10}{10}}
\put(2963,1737){\makebox(0,0)[lb]{\smash{{{\SetFigFont{5}{6.0}{\rmdefault}{\mddefault}{\updefault}$0$}}}}}
\put(2812,1437){\makebox(0,0)[lb]{\smash{{{\SetFigFont{5}{6.0}{\rmdefault}{\mddefault}{\updefault}$0$}}}}}
\put(2737,1287){\makebox(0,0)[lb]{\smash{{{\SetFigFont{5}{6.0}{\rmdefault}{\mddefault}{\updefault}$1$}}}}}
\put(2587,987){\makebox(0,0)[lb]{\smash{{{\SetFigFont{5}{6.0}{\rmdefault}{\mddefault}{\updefault}$1$}}}}}
\put(2250,312){\blacken\ellipse{10}{10}}
\put(2250,312){\ellipse{10}{10}}
\put(2212,237){\blacken\ellipse{10}{10}}
\put(2212,237){\ellipse{10}{10}}
\put(2175,162){\blacken\ellipse{10}{10}}
\put(2175,162){\ellipse{10}{10}}
\put(2475,762){\blacken\ellipse{10}{10}}
\put(2475,762){\ellipse{10}{10}}
\put(2437,687){\blacken\ellipse{10}{10}}
\put(2437,687){\ellipse{10}{10}}
\put(2400,612){\blacken\ellipse{10}{10}}
\put(2400,612){\ellipse{10}{10}}
\put(2512,837){\makebox(0,0)[lb]{\smash{{{\SetFigFont{5}{6.0}{\rmdefault}{\mddefault}{\updefault}$0$}}}}}
\put(2362,537){\makebox(0,0)[lb]{\smash{{{\SetFigFont{5}{6.0}{\rmdefault}{\mddefault}{\updefault}$0$}}}}}
\put(2287,387){\makebox(0,0)[lb]{\smash{{{\SetFigFont{5}{6.0}{\rmdefault}{\mddefault}{\updefault}$1$}}}}}
\put(2137,87){\makebox(0,0)[lb]{\smash{{{\SetFigFont{5}{6.0}{\rmdefault}{\mddefault}{\updefault}$1$}}}}}
\put(3300,1811){\makebox(0,0)[lb]{\smash{{{\SetFigFont{8}{9.6}{\rmdefault}{\mddefault}{\updefault}with pieces}}}}}
\put(3300,1661){\makebox(0,0)[lb]{\smash{{{\SetFigFont{8}{9.6}{\rmdefault}{\mddefault}{\updefault}of the form}}}}}
\put(0,1886){\makebox(0,0)[lb]{\smash{{{\SetFigFont{8}{9.6}{\rmdefault}{\mddefault}{\updefault}ways can you}}}}}
\put(0,1736){\makebox(0,0)[lb]{\smash{{{\SetFigFont{8}{9.6}{\rmdefault}{\mddefault}{\updefault}fill in the}}}}}
\put(0,1586){\makebox(0,0)[lb]{\smash{{{\SetFigFont{8}{9.6}{\rmdefault}{\mddefault}{\updefault}puzzle}}}}}
\put(0,2036){\makebox(0,0)[lb]{\smash{{{\SetFigFont{8}{9.6}{\rmdefault}{\mddefault}{\updefault}In how many}}}}}
\path(1162,1811)(2062,3611)(2963,1811)
	(2062,12)(1162,1811)(2963,1811)
\path(4313,1587)(4388,1737)(4463,1587)(4313,1587)
\path(4763,1886)(4838,2036)(4913,1886)
	(4838,1737)(4763,1886)
\path(4313,1961)(4388,2111)(4463,1961)(4313,1961)
\put(5325,1811){\makebox(0,0)[lb]{\smash{{{\SetFigFont{8}{9.6}{\rmdefault}{\mddefault}{\updefault}?}}}}}
\put(4913,1961){\makebox(0,0)[lb]{\smash{{{\SetFigFont{5}{6.0}{\rmdefault}{\mddefault}{\updefault}$0$}}}}}
\put(4440,1662){\makebox(0,0)[lb]{\smash{{{\SetFigFont{5}{6.0}{\rmdefault}{\mddefault}{\updefault}$1$}}}}}
\put(4440,2036){\makebox(0,0)[lb]{\smash{{{\SetFigFont{5}{6.0}{\rmdefault}{\mddefault}{\updefault}$0$}}}}}
\put(4916,1774){\makebox(0,0)[lb]{\smash{{{\SetFigFont{5}{6.0}{\rmdefault}{\mddefault}{\updefault}$1$}}}}}
\put(4363,1881){\makebox(0,0)[lb]{\smash{{{\SetFigFont{5}{6.0}{\rmdefault}{\mddefault}{\updefault}$0$}}}}}
\put(4359,1498){\makebox(0,0)[lb]{\smash{{{\SetFigFont{5}{6.0}{\rmdefault}{\mddefault}{\updefault}$1$}}}}}
\put(4708,1785){\makebox(0,0)[lb]{\smash{{{\SetFigFont{5}{6.0}{\rmdefault}{\mddefault}{\updefault}$0$}}}}}
\put(4702,1962){\makebox(0,0)[lb]{\smash{{{\SetFigFont{5}{6.0}{\rmdefault}{\mddefault}{\updefault}$1$}}}}}
\put(4293,1665){\makebox(0,0)[lb]{\smash{{{\SetFigFont{5}{6.0}{\rmdefault}{\mddefault}{\updefault}$1$}}}}}
\put(4286,2037){\makebox(0,0)[lb]{\smash{{{\SetFigFont{5}{6.0}{\rmdefault}{\mddefault}{\updefault}$0$}}}}}
\end{picture}
}

%% file: F.tex
\setlength{\unitlength}{0.00083333in}
\begingroup\makeatletter\ifx\SetFigFont\undefined%
\gdef\SetFigFont#1#2#3#4#5{%
  \reset@font\fontsize{#1}{#2pt}%
  \fontfamily{#3}\fontseries{#4}\fontshape{#5}%
  \selectfont}%
\fi\endgroup%
{\renewcommand{\dashlinestretch}{30}
\begin{picture}(1897,2043)(0,-10)
\put(1621,406){\blacken\ellipse{10}{10}}
\put(1621,406){\ellipse{10}{10}}
\put(1584,331){\blacken\ellipse{10}{10}}
\put(1584,331){\ellipse{10}{10}}
\put(1546,256){\blacken\ellipse{10}{10}}
\put(1546,256){\ellipse{10}{10}}
\put(1846,856){\blacken\ellipse{10}{10}}
\put(1846,856){\ellipse{10}{10}}
\put(1809,781){\blacken\ellipse{10}{10}}
\put(1809,781){\ellipse{10}{10}}
\put(1771,706){\blacken\ellipse{10}{10}}
\put(1771,706){\ellipse{10}{10}}
\put(1884,931){\makebox(0,0)[lb]{\smash{{{\SetFigFont{5}{6.0}{\rmdefault}{\mddefault}{\updefault}$0$}}}}}
\put(1734,631){\makebox(0,0)[lb]{\smash{{{\SetFigFont{5}{6.0}{\rmdefault}{\mddefault}{\updefault}$0$}}}}}
\put(1659,481){\makebox(0,0)[lb]{\smash{{{\SetFigFont{5}{6.0}{\rmdefault}{\mddefault}{\updefault}$1$}}}}}
\put(1509,181){\makebox(0,0)[lb]{\smash{{{\SetFigFont{5}{6.0}{\rmdefault}{\mddefault}{\updefault}$1$}}}}}
\put(1246,1944){\blacken\ellipse{10}{10}}
\put(1246,1944){\ellipse{10}{10}}
\put(1171,1944){\blacken\ellipse{10}{10}}
\put(1171,1944){\ellipse{10}{10}}
\put(1096,1944){\blacken\ellipse{10}{10}}
\put(1096,1944){\ellipse{10}{10}}
\put(647,1944){\blacken\ellipse{10}{10}}
\put(647,1944){\ellipse{10}{10}}
\put(797,1944){\blacken\ellipse{10}{10}}
\put(797,1944){\ellipse{10}{10}}
\put(722,1944){\blacken\ellipse{10}{10}}
\put(722,1944){\ellipse{10}{10}}
\put(1021,1944){\makebox(0,0)[lb]{\smash{{{\SetFigFont{5}{6.0}{\rmdefault}{\mddefault}{\updefault}$1$}}}}}
\put(572,1944){\makebox(0,0)[lb]{\smash{{{\SetFigFont{5}{6.0}{\rmdefault}{\mddefault}{\updefault}$0$}}}}}
\put(872,1944){\makebox(0,0)[lb]{\smash{{{\SetFigFont{5}{6.0}{\rmdefault}{\mddefault}{\updefault}$0$}}}}}
\put(1321,1944){\makebox(0,0)[lb]{\smash{{{\SetFigFont{5}{6.0}{\rmdefault}{\mddefault}{\updefault}$1$}}}}}
\put(647,69){\blacken\ellipse{10}{10}}
\put(647,69){\ellipse{10}{10}}
\put(797,69){\blacken\ellipse{10}{10}}
\put(797,69){\ellipse{10}{10}}
\put(722,69){\blacken\ellipse{10}{10}}
\put(722,69){\ellipse{10}{10}}
\put(1246,69){\blacken\ellipse{10}{10}}
\put(1246,69){\ellipse{10}{10}}
\put(1171,69){\blacken\ellipse{10}{10}}
\put(1171,69){\ellipse{10}{10}}
\put(1096,69){\blacken\ellipse{10}{10}}
\put(1096,69){\ellipse{10}{10}}
\put(1584,1756){\blacken\ellipse{10}{10}}
\put(1584,1756){\ellipse{10}{10}}
\put(1621,1681){\blacken\ellipse{10}{10}}
\put(1621,1681){\ellipse{10}{10}}
\put(1659,1606){\blacken\ellipse{10}{10}}
\put(1659,1606){\ellipse{10}{10}}
\put(1846,1231){\blacken\ellipse{10}{10}}
\put(1846,1231){\ellipse{10}{10}}
\put(1884,1156){\blacken\ellipse{10}{10}}
\put(1884,1156){\ellipse{10}{10}}
\put(1809,1306){\blacken\ellipse{10}{10}}
\put(1809,1306){\ellipse{10}{10}}
\put(1509,1831){\makebox(0,0)[lb]{\smash{{{\SetFigFont{5}{6.0}{\rmdefault}{\mddefault}{\updefault}$0$}}}}}
\put(1884,1081){\makebox(0,0)[lb]{\smash{{{\SetFigFont{5}{6.0}{\rmdefault}{\mddefault}{\updefault}$1$}}}}}
\put(1734,1381){\makebox(0,0)[lb]{\smash{{{\SetFigFont{5}{6.0}{\rmdefault}{\mddefault}{\updefault}$1$}}}}}
\put(1659,1531){\makebox(0,0)[lb]{\smash{{{\SetFigFont{5}{6.0}{\rmdefault}{\mddefault}{\updefault}$0$}}}}}
\path(94,1057)(994,1057)
\path(94,1097)(94,1019)
\path(994,1095)(994,1019)
\put(84,856){\blacken\ellipse{10}{10}}
\put(84,856){\ellipse{10}{10}}
\put(122,781){\blacken\ellipse{10}{10}}
\put(122,781){\ellipse{10}{10}}
\put(159,706){\blacken\ellipse{10}{10}}
\put(159,706){\ellipse{10}{10}}
\put(309,406){\blacken\ellipse{10}{10}}
\put(309,406){\ellipse{10}{10}}
\put(347,331){\blacken\ellipse{10}{10}}
\put(347,331){\ellipse{10}{10}}
\put(384,256){\blacken\ellipse{10}{10}}
\put(384,256){\ellipse{10}{10}}
\put(47,931){\makebox(0,0)[lb]{\smash{{{\SetFigFont{5}{6.0}{\rmdefault}{\mddefault}{\updefault}$1$}}}}}
\put(272,481){\makebox(0,0)[lb]{\smash{{{\SetFigFont{5}{6.0}{\rmdefault}{\mddefault}{\updefault}$0$}}}}}
\put(422,181){\makebox(0,0)[lb]{\smash{{{\SetFigFont{5}{6.0}{\rmdefault}{\mddefault}{\updefault}$0$}}}}}
\put(197,631){\makebox(0,0)[lb]{\smash{{{\SetFigFont{5}{6.0}{\rmdefault}{\mddefault}{\updefault}$1$}}}}}
\put(384,1756){\blacken\ellipse{10}{10}}
\put(384,1756){\ellipse{10}{10}}
\put(347,1681){\blacken\ellipse{10}{10}}
\put(347,1681){\ellipse{10}{10}}
\put(309,1606){\blacken\ellipse{10}{10}}
\put(309,1606){\ellipse{10}{10}}
\put(159,1306){\blacken\ellipse{10}{10}}
\put(159,1306){\ellipse{10}{10}}
\put(122,1231){\blacken\ellipse{10}{10}}
\put(122,1231){\ellipse{10}{10}}
\put(84,1156){\blacken\ellipse{10}{10}}
\put(84,1156){\ellipse{10}{10}}
\path(84,1006)(534,1906)(1434,1906)
	(1884,1006)(1434,106)(534,106)
	(84,1006)(1884,1006)
\put(422,1831){\makebox(0,0)[lb]{\smash{{{\SetFigFont{5}{6.0}{\rmdefault}{\mddefault}{\updefault}$1$}}}}}
\put(272,1531){\makebox(0,0)[lb]{\smash{{{\SetFigFont{5}{6.0}{\rmdefault}{\mddefault}{\updefault}$1$}}}}}
\put(197,1381){\makebox(0,0)[lb]{\smash{{{\SetFigFont{5}{6.0}{\rmdefault}{\mddefault}{\updefault}$0$}}}}}
\put(572,31){\makebox(0,0)[lb]{\smash{{{\SetFigFont{5}{6.0}{\rmdefault}{\mddefault}{\updefault}$1$}}}}}
\put(1021,31){\makebox(0,0)[lb]{\smash{{{\SetFigFont{5}{6.0}{\rmdefault}{\mddefault}{\updefault}$0$}}}}}
\put(872,31){\makebox(0,0)[lb]{\smash{{{\SetFigFont{5}{6.0}{\rmdefault}{\mddefault}{\updefault}$1$}}}}}
\put(1321,31){\makebox(0,0)[lb]{\smash{{{\SetFigFont{5}{6.0}{\rmdefault}{\mddefault}{\updefault}$0$}}}}}
\put(589,1094){\makebox(0,0)[lb]{\smash{{{\SetFigFont{5}{6.0}{\rmdefault}{\mddefault}{\updefault}$A$}}}}}
\put(0,1035){\makebox(0,0)[lb]{\smash{{{\SetFigFont{5}{6.0}{\rmdefault}{\mddefault}{\updefault}$0$}}}}}
\end{picture}
}

%% file: G.tex
\setlength{\unitlength}{0.00083333in}
\begingroup\makeatletter\ifx\SetFigFont\undefined%
\gdef\SetFigFont#1#2#3#4#5{%
  \reset@font\fontsize{#1}{#2pt}%
  \fontfamily{#3}\fontseries{#4}\fontshape{#5}%
  \selectfont}%
\fi\endgroup%
{\renewcommand{\dashlinestretch}{30}
\begin{picture}(1412,831)(0,-10)
\path(12,612)(87,612)
\path(50,612)(50,537)
\path(12,312)(87,312)
\path(50,387)(50,312)
\path(725,762)(800,762)
\path(200,762)(275,762)
\path(200,800)(200,725)
\path(800,800)(800,725)
\path(1325,762)(1400,762)
\path(800,762)(875,762)
\path(800,800)(800,725)
\path(1400,800)(1400,725)
\path(12,312)(87,312)
\path(50,312)(50,237)
\path(12,12)(87,12)
\path(50,87)(50,12)
\path(200,612)(1400,612)(1400,312)
	(800,312)(800,12)(200,12)
	(200,612)(200,612)
\path(200,312)(800,312)(800,612)
\put(30,425){\makebox(0,0)[lb]{\smash{{{\SetFigFont{5}{6.0}{\rmdefault}{\mddefault}{\updefault}$K$}}}}}
\put(350,732){\makebox(0,0)[lb]{\smash{{{\SetFigFont{5}{6.0}{\rmdefault}{\mddefault}{\updefault}$N-K$}}}}}
\put(950,732){\makebox(0,0)[lb]{\smash{{{\SetFigFont{5}{6.0}{\rmdefault}{\mddefault}{\updefault}$N-K$}}}}}
\put(30,125){\makebox(0,0)[lb]{\smash{{{\SetFigFont{5}{6.0}{\rmdefault}{\mddefault}{\updefault}$K$}}}}}
\put(1400,612){\makebox(0,0)[lb]{\smash{{{\SetFigFont{5}{6.0}{\rmdefault}{\mddefault}{\updefault}$3$}}}}}
\end{picture}
}

%% file: H.tex
\setlength{\unitlength}{0.00083333in}
\begingroup\makeatletter\ifx\SetFigFont\undefined%
\gdef\SetFigFont#1#2#3#4#5{%
  \reset@font\fontsize{#1}{#2pt}%
  \fontfamily{#3}\fontseries{#4}\fontshape{#5}%
  \selectfont}%
\fi\endgroup%
{\renewcommand{\dashlinestretch}{30}
\begin{picture}(2726,2813)(0,-10)
\put(1957,1740){\blacken\ellipse{10}{10}}
\put(1957,1740){\ellipse{10}{10}}
\put(1993,1665){\blacken\ellipse{10}{10}}
\put(1993,1665){\ellipse{10}{10}}
\put(2030,1590){\blacken\ellipse{10}{10}}
\put(2030,1590){\ellipse{10}{10}}
\put(2215,1216){\blacken\ellipse{10}{10}}
\put(2215,1216){\ellipse{10}{10}}
\put(2251,1141){\blacken\ellipse{10}{10}}
\put(2251,1141){\ellipse{10}{10}}
\put(2178,1290){\blacken\ellipse{10}{10}}
\put(2178,1290){\ellipse{10}{10}}
\put(1883,1814){\makebox(0,0)[lb]{\smash{{{\SetFigFont{5}{6.0}{\rmdefault}{\mddefault}{\updefault}$0$}}}}}
\put(2251,1066){\makebox(0,0)[lb]{\smash{{{\SetFigFont{5}{6.0}{\rmdefault}{\mddefault}{\updefault}$1$}}}}}
\put(2104,1365){\makebox(0,0)[lb]{\smash{{{\SetFigFont{5}{6.0}{\rmdefault}{\mddefault}{\updefault}$1$}}}}}
\put(2030,1515){\makebox(0,0)[lb]{\smash{{{\SetFigFont{5}{6.0}{\rmdefault}{\mddefault}{\updefault}$0$}}}}}
\put(2419,834){\blacken\ellipse{10}{10}}
\put(2419,834){\ellipse{10}{10}}
\put(2455,759){\blacken\ellipse{10}{10}}
\put(2455,759){\ellipse{10}{10}}
\put(2491,684){\blacken\ellipse{10}{10}}
\put(2491,684){\ellipse{10}{10}}
\put(2676,311){\blacken\ellipse{10}{10}}
\put(2676,311){\ellipse{10}{10}}
\put(2713,236){\blacken\ellipse{10}{10}}
\put(2713,236){\ellipse{10}{10}}
\put(2639,385){\blacken\ellipse{10}{10}}
\put(2639,385){\ellipse{10}{10}}
\put(2345,909){\makebox(0,0)[lb]{\smash{{{\SetFigFont{5}{6.0}{\rmdefault}{\mddefault}{\updefault}$0$}}}}}
\put(2713,161){\makebox(0,0)[lb]{\smash{{{\SetFigFont{5}{6.0}{\rmdefault}{\mddefault}{\updefault}$1$}}}}}
\put(2565,460){\makebox(0,0)[lb]{\smash{{{\SetFigFont{5}{6.0}{\rmdefault}{\mddefault}{\updefault}$1$}}}}}
\put(2491,610){\makebox(0,0)[lb]{\smash{{{\SetFigFont{5}{6.0}{\rmdefault}{\mddefault}{\updefault}$0$}}}}}
\put(775,1740){\blacken\ellipse{10}{10}}
\put(775,1740){\ellipse{10}{10}}
\put(738,1665){\blacken\ellipse{10}{10}}
\put(738,1665){\ellipse{10}{10}}
\put(701,1590){\blacken\ellipse{10}{10}}
\put(701,1590){\ellipse{10}{10}}
\put(555,1290){\blacken\ellipse{10}{10}}
\put(555,1290){\ellipse{10}{10}}
\put(518,1216){\blacken\ellipse{10}{10}}
\put(518,1216){\ellipse{10}{10}}
\put(481,1141){\blacken\ellipse{10}{10}}
\put(481,1141){\ellipse{10}{10}}
\put(2014,57){\blacken\ellipse{10}{10}}
\put(2014,57){\ellipse{10}{10}}
\put(1940,57){\blacken\ellipse{10}{10}}
\put(1940,57){\ellipse{10}{10}}
\put(2089,57){\blacken\ellipse{10}{10}}
\put(2089,57){\ellipse{10}{10}}
\put(2531,57){\blacken\ellipse{10}{10}}
\put(2531,57){\ellipse{10}{10}}
\put(2457,57){\blacken\ellipse{10}{10}}
\put(2457,57){\ellipse{10}{10}}
\put(2383,57){\blacken\ellipse{10}{10}}
\put(2383,57){\ellipse{10}{10}}
\put(1126,57){\blacken\ellipse{10}{10}}
\put(1126,57){\ellipse{10}{10}}
\put(1052,57){\blacken\ellipse{10}{10}}
\put(1052,57){\ellipse{10}{10}}
\put(1200,57){\blacken\ellipse{10}{10}}
\put(1200,57){\ellipse{10}{10}}
\put(1642,57){\blacken\ellipse{10}{10}}
\put(1642,57){\ellipse{10}{10}}
\put(1568,57){\blacken\ellipse{10}{10}}
\put(1568,57){\ellipse{10}{10}}
\put(1495,57){\blacken\ellipse{10}{10}}
\put(1495,57){\ellipse{10}{10}}
\put(238,57){\blacken\ellipse{10}{10}}
\put(238,57){\ellipse{10}{10}}
\put(164,57){\blacken\ellipse{10}{10}}
\put(164,57){\ellipse{10}{10}}
\put(312,57){\blacken\ellipse{10}{10}}
\put(312,57){\ellipse{10}{10}}
\put(754,57){\blacken\ellipse{10}{10}}
\put(754,57){\ellipse{10}{10}}
\put(680,57){\blacken\ellipse{10}{10}}
\put(680,57){\ellipse{10}{10}}
\put(607,57){\blacken\ellipse{10}{10}}
\put(607,57){\ellipse{10}{10}}
\put(111,393){\blacken\ellipse{10}{10}}
\put(111,393){\ellipse{10}{10}}
\put(74,319){\blacken\ellipse{10}{10}}
\put(74,319){\ellipse{10}{10}}
\put(37,244){\blacken\ellipse{10}{10}}
\put(37,244){\ellipse{10}{10}}
\put(333,841){\blacken\ellipse{10}{10}}
\put(333,841){\ellipse{10}{10}}
\put(296,767){\blacken\ellipse{10}{10}}
\put(296,767){\ellipse{10}{10}}
\put(259,692){\blacken\ellipse{10}{10}}
\put(259,692){\ellipse{10}{10}}
\put(1514,2637){\blacken\ellipse{10}{10}}
\put(1514,2637){\ellipse{10}{10}}
\put(1550,2562){\blacken\ellipse{10}{10}}
\put(1550,2562){\ellipse{10}{10}}
\put(1587,2487){\blacken\ellipse{10}{10}}
\put(1587,2487){\ellipse{10}{10}}
\put(1772,2113){\blacken\ellipse{10}{10}}
\put(1772,2113){\ellipse{10}{10}}
\put(1809,2039){\blacken\ellipse{10}{10}}
\put(1809,2039){\ellipse{10}{10}}
\put(1735,2188){\blacken\ellipse{10}{10}}
\put(1735,2188){\ellipse{10}{10}}
\put(997,2188){\blacken\ellipse{10}{10}}
\put(997,2188){\ellipse{10}{10}}
\put(960,2113){\blacken\ellipse{10}{10}}
\put(960,2113){\ellipse{10}{10}}
\put(923,2039){\blacken\ellipse{10}{10}}
\put(923,2039){\ellipse{10}{10}}
\put(1219,2637){\blacken\ellipse{10}{10}}
\put(1219,2637){\ellipse{10}{10}}
\put(1182,2562){\blacken\ellipse{10}{10}}
\put(1182,2562){\ellipse{10}{10}}
\put(1145,2487){\blacken\ellipse{10}{10}}
\put(1145,2487){\ellipse{10}{10}}
\path(37,94)(1367,2786)(2696,94)
	(1809,94)(37,94)
\put(444,1066){\makebox(0,0)[lb]{\smash{{{\SetFigFont{5}{6.0}{\rmdefault}{\mddefault}{\updefault}$0$}}}}}
\put(812,1814){\makebox(0,0)[lb]{\smash{{{\SetFigFont{5}{6.0}{\rmdefault}{\mddefault}{\updefault}$1$}}}}}
\put(665,1515){\makebox(0,0)[lb]{\smash{{{\SetFigFont{5}{6.0}{\rmdefault}{\mddefault}{\updefault}$1$}}}}}
\put(592,1365){\makebox(0,0)[lb]{\smash{{{\SetFigFont{5}{6.0}{\rmdefault}{\mddefault}{\updefault}$0$}}}}}
\put(960,19){\makebox(0,0)[lb]{\smash{{{\SetFigFont{5}{6.0}{\rmdefault}{\mddefault}{\updefault}$1$}}}}}
\put(2292,19){\makebox(0,0)[lb]{\smash{{{\SetFigFont{5}{6.0}{\rmdefault}{\mddefault}{\updefault}$0$}}}}}
\put(72,19){\makebox(0,0)[lb]{\smash{{{\SetFigFont{5}{6.0}{\rmdefault}{\mddefault}{\updefault}$1$}}}}}
\put(810,19){\makebox(0,0)[lb]{\smash{{{\SetFigFont{5}{6.0}{\rmdefault}{\mddefault}{\updefault}$0$}}}}}
\put(1404,19){\makebox(0,0)[lb]{\smash{{{\SetFigFont{5}{6.0}{\rmdefault}{\mddefault}{\updefault}$0$}}}}}
\put(1256,19){\makebox(0,0)[lb]{\smash{{{\SetFigFont{5}{6.0}{\rmdefault}{\mddefault}{\updefault}$1$}}}}}
\put(515,19){\makebox(0,0)[lb]{\smash{{{\SetFigFont{5}{6.0}{\rmdefault}{\mddefault}{\updefault}$0$}}}}}
\put(368,19){\makebox(0,0)[lb]{\smash{{{\SetFigFont{5}{6.0}{\rmdefault}{\mddefault}{\updefault}$1$}}}}}
\put(1698,19){\makebox(0,0)[lb]{\smash{{{\SetFigFont{5}{6.0}{\rmdefault}{\mddefault}{\updefault}$0$}}}}}
\put(1848,19){\makebox(0,0)[lb]{\smash{{{\SetFigFont{5}{6.0}{\rmdefault}{\mddefault}{\updefault}$1$}}}}}
\put(2144,19){\makebox(0,0)[lb]{\smash{{{\SetFigFont{5}{6.0}{\rmdefault}{\mddefault}{\updefault}$1$}}}}}
\put(2586,19){\makebox(0,0)[lb]{\smash{{{\SetFigFont{5}{6.0}{\rmdefault}{\mddefault}{\updefault}$0$}}}}}
\put(0,169){\makebox(0,0)[lb]{\smash{{{\SetFigFont{5}{6.0}{\rmdefault}{\mddefault}{\updefault}$0$}}}}}
\put(148,468){\makebox(0,0)[lb]{\smash{{{\SetFigFont{5}{6.0}{\rmdefault}{\mddefault}{\updefault}$0$}}}}}
\put(370,916){\makebox(0,0)[lb]{\smash{{{\SetFigFont{5}{6.0}{\rmdefault}{\mddefault}{\updefault}$1$}}}}}
\put(222,618){\makebox(0,0)[lb]{\smash{{{\SetFigFont{5}{6.0}{\rmdefault}{\mddefault}{\updefault}$1$}}}}}
\put(1440,2711){\makebox(0,0)[lb]{\smash{{{\SetFigFont{5}{6.0}{\rmdefault}{\mddefault}{\updefault}$0$}}}}}
\put(1809,1964){\makebox(0,0)[lb]{\smash{{{\SetFigFont{5}{6.0}{\rmdefault}{\mddefault}{\updefault}$1$}}}}}
\put(1661,2262){\makebox(0,0)[lb]{\smash{{{\SetFigFont{5}{6.0}{\rmdefault}{\mddefault}{\updefault}$1$}}}}}
\put(1587,2412){\makebox(0,0)[lb]{\smash{{{\SetFigFont{5}{6.0}{\rmdefault}{\mddefault}{\updefault}$0$}}}}}
\put(886,1964){\makebox(0,0)[lb]{\smash{{{\SetFigFont{5}{6.0}{\rmdefault}{\mddefault}{\updefault}$0$}}}}}
\put(1034,2262){\makebox(0,0)[lb]{\smash{{{\SetFigFont{5}{6.0}{\rmdefault}{\mddefault}{\updefault}$0$}}}}}
\put(1108,2412){\makebox(0,0)[lb]{\smash{{{\SetFigFont{5}{6.0}{\rmdefault}{\mddefault}{\updefault}$1$}}}}}
\put(1256,2711){\makebox(0,0)[lb]{\smash{{{\SetFigFont{5}{6.0}{\rmdefault}{\mddefault}{\updefault}$1$}}}}}
\end{picture}
}